%% file: arxiv.tex
\documentclass[11pt]{article}

\usepackage[hyperfootnotes = true,colorlinks,citecolor=lightblue,urlcolor=lightblue, linkcolor=lightblue]{hyperref}

\usepackage{enumerate}
\usepackage{xcolor}
\definecolor{lightblue}{HTML}{044E9E}
\definecolor{cblue}{HTML}{FD840B}
\definecolor{newblue}{HTML}{164472}
\definecolor{light-gray}{gray}{0.4}
\definecolor{example-color}{HTML}{164472}

\usepackage{amsmath}
\usepackage{graphicx}
\usepackage{psfrag}
\usepackage{epsf}
\usepackage{enumerate}
\usepackage{enumitem}
\usepackage[round,longnamesfirst]{natbib}
\usepackage{url} 
\usepackage{amsthm}
\usepackage{xpatch}
\usepackage{subcaption}

\usepackage{longtable}
\usepackage{mathrsfs}
\usepackage{tikz}
\usepackage{color}
\usepackage{amssymb}
\usepackage{latexsym}
\usepackage{xr}
\usepackage{epstopdf}
\usepackage{float}
\usepackage{multirow}
\usepackage{dsfont}
\usepackage{mathabx}

\numberwithin{equation}{section}

\usepackage{sectsty}
\subsectionfont{\normalfont\itshape}
\subsubsectionfont{\normalfont\itshape}

\def\RR{\mathbb R}
\def\YY{\mathbb Y}
\def\ZZ{\mathbb Z}

\def\NN{\mathbb N}

\def\P{\operatorname P}

\newcommand{\EE}{\operatorname{E}} 
\newcommand{\Var}{\operatorname{Var}}

\newcommand{\op}{\operatorname{op}}

\newcommand{\distr}{\operatorname{d}} 

\newcommand{\Hi}{\mathbb{H}}

\newcommand{\esssup}{\operatornamewithlimits{ess\text{ }sup}}

\newtheorem{lemma}{Lemma}[section]

\newtheorem{theorem}[lemma]{Theorem}

\newtheorem{definition}[lemma]{Definition}
\newtheorem{corollary}[lemma]{Corollary}

\newtheorem{remark}[lemma]{Remark}

\newtheorem{openquest}{Question}

\makeatletter
\xpatchcmd{\proof}{\@addpunct{.}}{\@addpunct{:}}{}{}
\makeatother

\DeclareFontFamily{U}{mathx}{\hyphenchar\font45}
\DeclareFontShape{U}{mathx}{m}{n}{<-> mathx10}{}
\DeclareSymbolFont{mathx}{U}{mathx}{m}{n}
\DeclareMathAccent{\widebar}{0}{mathx}{"73}

\usepackage{color}
\usepackage{marginnote}

\usepackage{geometry}
\allowdisplaybreaks

\def\spacingset#1{\renewcommand{\baselinestretch}%
{#1}\small\normalsize} \spacingset{1}

\begin{document}
%

\title{{\bf Scaling limits for sample autocovariance operators of Hilbert space-valued linear processes}\footnote{An earlier report on this problem titled ``Sample autocovariance operators of long-range dependent Hilbert space-valued linear processes" had appeared on MD's website in 2018.}}

\author{
Marie-Christine D\"uker \\ FAU Erlangen-N\"urnberg                             \and
Pavlos Zoubouloglou      \\ University of M\"unster}

\maketitle


%
%

\maketitle

\begin{abstract}
\noindent
This article considers linear processes with values in a separable Hilbert space exhibiting long-range dependence. The scaling limits for the sample autocovariance operators at different time lags are investigated in the topology of their respective Hilbert spaces. Distinguishing two different regimes of long-range dependence, the limiting object is either a Hilbert space-valued Gaussian or a Hilbert space-valued non-Gaussian random variable. The latter can be represented as a unitary transformation of double Wiener-It\^o integrals with sample paths in a function space.
This work is the first to show weak convergence to such double stochastic integrals with sample paths in infinite dimensions. The result generalizes the well known convergence to a Hermite process in finite dimensions, introducing a new domain of attraction for probability measures in Hilbert spaces.
The key technical contributions include the introduction of double Wiener-It\^o integrals with values in a function space and with dependent integrators, as well as establishing sufficient conditions for their existence as limits of sample autocovariance operators.


\end{abstract}

\input{main}

\input{appendix}

%
%
%
%
%
%
%
%
%
%

\bigskip
\noindent {\bf Acknowledgments:} MD's research was supported by the the Deutsche Forschungsgemeinschaft (DFG, German Research Foundation) - Research Training Group 2131 - \textit{High-dimensional Phenomena in Probability - Fluctuations and Discontinuity}. PZ's research was funded by the DFG under Germany's Excellence Strategy EXC 2044 –390685587, Mathematics Münster: Dynamics–Geometry–Structure. 
\bigskip

\bibliographystyle{plainnat}
\bibliography{sampleautocovhilbert}

\newpage

\end{document}

%% file: main.tex
\section{Introduction}
\label{s:introduction}
In this article, we investigate the weak convergence of the sample autocovariance operators of a Hilbert space-valued linear process exhibiting long-range dependence. Our setting is as follows: Let $\Hi$ denote a separable Hilbert space equipped with the inner product $\langle \cdot , \cdot \rangle_\Hi$ and norm $\| \cdot \|_\Hi$. We further write $L(\Hi)$ for the set of all bounded linear operators on $\Hi$.
We consider a sequence of random variables $\{X_{n}\}_{n \in \ZZ}$ defined on some probability space $(\Omega,\mathcal{F},\P)$ with values in $\Hi$. 
Suppose the stochastic process $\{X_{n}\}_{n \in \ZZ}$ admits the linear representation
\begin{align} \label{eq:linear_process}
X_{n}=\sum_{j = 0}^{\infty} u_{j}[\varepsilon_{n-j}], \hspace{0.2cm} n \in \ZZ,
\end{align}
with $u_{j} \in L(\Hi)$ for all $j \in \NN$ and $\{ \varepsilon_{j} \}_{j\in\ZZ}$ is a sequence of $\Hi$-valued independent, identically distributed (i.i.d.), zero-mean random variables.

For a given time lag $h \in \NN_0$, the sample autocovariance operator $\widehat{\Gamma}_{N,h}$ and its population counterpart $\Gamma_h$ of a process $\{X_n\}_{n \in \ZZ}$ are given by
\begin{equation} \label{eq:SampleAutoShort}
\widehat{\Gamma}_{N,h} \doteq \frac{1}{N} \sum_{n=1}^{N} X_{n+h} \otimes X_{n} 
\hspace{0.2cm}
\text{ and }
\hspace{0.2cm}
\Gamma_{h} \doteq \EE ( X_{h} \otimes X_{0}).
\end{equation}
Both quantities in \eqref{eq:SampleAutoShort} take values in $\Hi^{\otimes 2}$.
We are interested in the weak convergence of the normalized operators 
\begin{equation} \label{eq:sampleautocvo2}
   A_N((\widehat{\Gamma}_{N,h}-\Gamma_{h}), h=0,\dots,H),  \quad \text{as }N \to \infty,
\end{equation}
in $(\Hi^{\otimes 2})^{\times (H+1)}$, where $\{A_N\}_{N \in \NN} \in L((\Hi^{\otimes 2})^{\times (H+1)})$ is a sequence of suitable scaling operators.

The asymptotic behavior of \eqref{eq:SampleAutoShort}--\eqref{eq:sampleautocvo2} is by now well understood when $\{X_n\}$ in \eqref{eq:linear_process} (for general $\{u_j\}_{j \in \NN_0}$) exhibits Short-Range Dependence (SRD), i.e., when $\sum_{j=0}^\infty \| u_j \|_{\operatorname{op}} < \infty$, where $\|\cdot \|_{\operatorname{op}}$ denotes the usual operator norm recalled in \eqref{eq:def-op-norm} below. In this case, the normalizing sequence is simply $A_N = N^{1/2}$ and the quantity in \eqref{eq:sampleautocvo2} satisfies a Central Limit Theorem (CLT) in $\Hi \otimes \Hi$ (or, equivalently, the space of Hilbert-Schmidt operators on $\Hi$), i.e., converges to a centered Gaussian random variable taking values in $\Hi \otimes \Hi$, with an explicit covariance operator; see \cite{Mas} and \eqref{eq:covar-oper-Sigma}--\eqref{def:def:Aq-C} below.

In contrast to \cite{Mas}, we are interested in the case where $\{X_n\}$ in \eqref{eq:linear_process} exhibits Long-Range Dependence (LRD), in the sense that,
\begin{equation} \label{eq:op-lrd}
    \sum_{j=0}^\infty \| u_j \|_{\operatorname{op}} = \infty.
\end{equation}
Under \eqref{eq:op-lrd}, the fluctuations of the sample autocovariance operator \eqref{eq:sampleautocvo2} crucially depend on the convergence or divergence of the series $\sum_{j=0}^{\infty} \| u_j \|^{4/3}_{\operatorname{op}}$. When this series converges, we obtain a CLT for the sample autocovariance operators for general $\{u_j\}_{j \in \NN_0}$. When it diverges, in addition to \eqref{eq:linear_process}, we need to impose the additional structure on $\{u_j\}_{j \in \NN_0}$ given by
\begin{align} \label{eq:linear_process2}
u_j \doteq (j+1)^{T-I}, \; j \in \NN, 
\end{align}
where $I$ is the identity operator in $\Hi$ and $T \in L(\Hi)$ is a self-adjoint operator; see Section \ref{subsec:prelim-linear} below for more details. In this case, both the usual scaling $A_N = N^{1/2}$ and the Gaussian limit are lost.

Unless otherwise stated and to simplify the presentation, we assume for the rest of the introduction that $\Hi$ is the space of square integrable functions with regard to a $\sigma$-finite measure $\mu$, i.e., $\mathbb{H} = L^2(\mathbb{Y}, \mathcal{A}, \mu)$. We assume, moreover, that $T = D_d$ is a multiplication operator associated to a measurable function $d$, evaluated by $D_{d}f \doteq \{ d(y) f(y), y \in \mathbb{Y}\}$ for all $f \in L^2(\mathbb{Y}, \mathcal{A}, \mu)$. Then, \eqref{eq:linear_process2} is recast as
\begin{equation} \label{eq:uj-lrd}
    u_j \doteq (j+1)^{D_d- I}.
\end{equation}
When $d(s) \in \left(0, \frac{1}{2} \right)$ for $s \in \mathbb{Y}$, $\{ u_j \}$ in \eqref{eq:uj-lrd} satisfies \eqref{eq:op-lrd}, which justifies referring to \eqref{eq:linear_process} with \eqref{eq:uj-lrd} and $d(s) \in \left(0, \frac{1}{2} \right)$ as the LRD case. Therefore, the results of \cite{Mas}, that require absolute summability of $\{u_j\}$ in the operator norm, are no longer applicable.

The sample mean for LRD case \eqref{eq:linear_process} with \eqref{eq:uj-lrd} (i.e., $\mathbb{H} = L^2(\mathbb{Y}, \mathcal{A}, \mu)$) was studied in \cite{Rackauskas2010,Rack_Suquet_Op_2011}. The authors derived a Gaussian limiting distribution for the sample mean and the piecewise linear functions associated with the partial sums. Their results were generalized by \citet{DUEKER} to LRD processes with values in a general Hilbert space $\Hi$ and representation \eqref{eq:linear_process} with \eqref{eq:linear_process2}. The fluctuations of the sample autocovariance operator for this model have not been studied, and are the main objective of the present paper.

LRD has been studied extensively for stochastic processes in finite dimensions; see \citet{giraitis,beran2013long,PipirasTaqqu} to name a few. Special to long-range dependence is the Rosenblatt process, a non-Gaussian process that can be represented as a double Wiener-It\^o  integral and arises as a limiting object for stationary processes exhibiting long-range dependence; see \citet{Taqqu1975,Rosenblatt1979} for some of the earlier references and \cite{tudor,veillette2013,leonenko2017,bai2017} to name a few more recent ones.

The weak convergence of the (suitably re-normalized) sample autocovariance estimators in \eqref{eq:SampleAutoShort}--\eqref{eq:sampleautocvo2} when the underlying model \eqref{eq:linear_process} with \eqref{eq:uj-lrd} takes values in $\Hi = \RR$ (and hence has a constant memory parameter $d$) was studied in \citet{HorvathKokoszka}. There, it was shown that the scaling limit is a Brownian motion with $A_N = N^{1/2}$ when $d \in \left(0,\frac{1}{4}\right)$, and follows the Rosenblatt distribution with $A_N = N^{1-2d}$ when $d \in \left(\frac{1}{4},\frac{1}{2}\right)$. In \citet{Duker:multiMix2020}, these results were extended to multivariate linear processes, allowing some components to exhibit  short- and others long-range dependence. 

In analogy to \citet{HorvathKokoszka}, the two regimes $d \in \left(0,\frac{1}{4}\right)$ and $d \in \left(\frac{1}{4},\frac{1}{2}\right)$ for the fluctuations of the autocovariance operator roughly correspond to the two regimes $d(s) \in \left(0,\frac{1}{4}\right)$ and $d(s) \in \left(\frac{1}{4},\frac{1}{2}\right)$ for each $s \in \mathbb{Y}$ in our case. We prove weak convergence of \eqref{eq:SampleAutoShort} to a Gaussian law in $L^2(\mathbb{Y}, \mathcal{A},\mu)$ for the first regime, and to a double Wiener-It\^o  integral with sample paths in $L^2(\mathbb{Y}, \mathcal{A},\mu)$ for the second regime. For the first regime, we extend the arguments of \cite{Mas}. The second regime requires a careful analysis of double Wiener-It\^o  integrals with sample paths in Hilbert spaces and with spatially dependent integrators. This substantially extends results from \cite{fox1985} and \cite{Nor94}; see Section \ref{se:biW}. Establishing the tools that prove convergence to a non-Gaussian limit in functional spaces is one of the key challenges to be addressed in this work; see Lemma \ref{Lemma_off_diagonals}. To the best of our knowledge, this is the first instance of convergence to a double Wiener-It\^o integral with sample paths in a Hilbert space. 
\par
Our work is closely related to functional data analysis, where time series in infinite dimensions exhibiting LRD have appeared in numerous domains such as finance; see, e.g., \cite{Caj05,Alvrod08,Preciado2008,CasGao08}. It is therefore reasonable to model such problems using $\{X_n\}_{n \in \ZZ}$ as in \eqref{eq:linear_process} with \eqref{eq:linear_process2}. Limit theorems such as those developed in this work are important for applications, as they can be used to design hypothesis tests and construct confidence bands.
We emphasize that, although there are many works dealing with modeling and theoretical aspects of functional time series exhibiting SRD (e.g., \cite{Merlevede1997,jirak2018,Mas,DueZou24,RadKrePap24}), such works are comparatively scarce in the context of LRD. Notable exceptions include the works of \cite{Chara_Rack_Central,Chara_Rack_Operator,DUEKER,LiRobSha20,ruiz22,DurRou24}. For further details on functional data analysis, we refer to \cite{HsiEub15,hormann2018}, while \cite{Bosq2000:Linear} remains the canonical reference for linear processes with values in Banach spaces.
\par
The rest of the paper is organized as follows. In Section \ref{s:preliminaries}, we introduce some notation and recall some preliminary technical results. 
In Section \ref{s:mainresults}, we present our main results on the scaling limits of the sample autocovariance operators of the process $\{X_{n}\}_{n \in \ZZ }$. In Section \ref{se:biW}, we introduce double Wiener-It\^o integrals with values in a function space and spatially dependent integrators; such objects are crucial in representing the resulting limit in the second regime.
In Section \ref{s:proofs}, we present the proofs of our main results. Section \ref{se:CH3appendix} is concerned with some technical lemmas and their proofs.

\section{Preliminaries}
\label{s:preliminaries}

In this section, we introduce notation, collect some preliminary facts on operators and function spaces used throughout the paper, 
and give some properties of the linear process \eqref{eq:linear_process}.

\subsection{Notation and terminology} 

Let $(X, \mathcal{G}, \nu)$ be a $\sigma$-finite measure space and $(Y, \| \cdot \|_{Y})$ a normed space. Then $L^2(X:Y) \doteq L^2(X,\mathcal{G},\nu:Y)$ denotes the space of square integrable, measurable functions on $(X,\mathcal{G},\nu)$, taking values in $Y$. When $(Y, \| \cdot \|_{Y}) = (\RR, |\cdot|_{\RR})$, we write $L^2(X)$ or $L^2(X, \nu)$. Recall that $L^2(X)$ is equipped with the inner product 
\begin{equation*}
\langle f,g \rangle_{L^2(X)} = \int_{X} f(s) g(s) \nu(ds), \hspace*{0.2cm} f,g \in L^2(X),
\end{equation*}
and its induced norm $\| \cdot \|_{L^2(X)}$. Moreover, denote, for a fixed integer $m \geq 1$, the space $L^2(X^{\times m}) \doteq L^2(X^{\times m}, \mathcal{G}^{\otimes m}, \nu^{\otimes m})$ of square integrable real-valued functions in the product measure space $(X^{\times m}, \mathcal{G}^{\otimes m}, \nu^{\otimes m})$. 

For $a,b > 0$, define the beta function
\begin{equation}
    \operatorname{B}(a,b) \doteq \int_{0}^{1}x^{a-1}(1-x)^{b-1}dx =\int_{0}^{\infty}x^{a-1}(x+1)^{-(a+b)}dx, \quad a, b > 0.
\end{equation}
Moreover, for a given function $d: \YY \to \RR_+$, we define two functions $c: \YY \times \YY \to \RR$ and $c:\YY \to \RR$ by
\begin{equation} \label{equality_beta}
    c(r,s)
\doteq \operatorname{B}(d(r),d(s))
=\int_{0}^{\infty}x^{d(r)-1}(x+1)^{d(s)-1}dx, \quad \text{and} \quad c(r) \doteq c(r,r),
\end{equation} 
respectively. 

Throughout this work, we write $(\RR, \mathcal{B}, \lambda)$ to denote the real numbers, equipped with the Lebesgue measure $\lambda$, and the Borel sets (in the usual topology) $\mathcal{B}$.

\subsection{Elements of operator theory}

Let $L(X,Y)$ be the space of bounded linear operators from $X$ to $Y$, and let $L(X) \doteq L(X,X)$. The space $(L(\Hi), \|\cdot \|_{\operatorname{op}})$ forms a Banach space, with the operator norm defined by 
\begin{equation} \label{eq:def-op-norm}
    \|T\|_{\operatorname{op}} \doteq \inf \{c \ge 0: \|Tv \|_{\Hi} \le c \|v\|_{\Hi} \;\text{for all } v \in \Hi \} ,\; T \in L(\Hi).
\end{equation}

The sample autocovariance operators are considered to be random elements with values in the space of Hilbert-Schmidt operators on $\Hi$, denoted by $\operatorname{HS}(\Hi)$. 
A Hilbert-Schmidt operator $A:\Hi \to \Hi$ is a bounded operator with finite Hilbert-Schmidt norm
\begin{equation} \label{def:HS_norm}
\| A \|^2_{\operatorname{HS}(\Hi)} \doteq 
\sum_{i=1}^{\infty} \| A e_{i} \|^2_{\Hi}<\infty,
\end{equation}
where $ \{e_{i}\}_{i \in \NN} $ is an orthonormal basis of $\Hi$.
The space $\operatorname{HS}(\Hi)$ equipped with the inner product $\langle A, B \rangle_{\operatorname{HS}(\Hi)}= \sum_{i=1}^{\infty} \langle Ae_{i}, Be_{i} \rangle_{\Hi} $ and its induced norm $\| A \|_{\operatorname{HS}(\Hi)}$ is a separable Hilbert space itself. Recall the (isometric) isomorphism $\Hi \otimes \Hi \cong \operatorname{HS}(\Hi)$; see, e.g., \cite{MuaFukSri17}, pp. 32-33. Invoking this isomorphism, we either show convergence of the sample autocovariance operators in $\operatorname{HS}(\Hi)$ (Theorem \ref{thm1:srd} and Corollary \ref{cor:srd-L^2}) or $\Hi \otimes \Hi$ (Theorems \ref{le:limit_theorem_L2} and \ref{th:main_result}), but think of weak convergence as equivalent in $\Hi \otimes \Hi$ and $\operatorname{HS}(\Hi)$. 

Closely related is the Banach space of trace class operators, denoted by $\operatorname{Tr}(\Hi)$, and equipped with the norm 
\begin{equation} \label{def:trace_norm}
     \| T \|_{\operatorname{Tr}(\Hi)} = \sum_{i=1}^{\infty} \langle |T| u_i, u_i \rangle_{\Hi},
\end{equation}
with $|T| = \sqrt{T^{*}T}$, where $T^*$ denotes the adjoint of the operator $T$. Whenever the space $\Hi$ is easily inferred from context, we may write $\operatorname{Tr}$ instead of $\operatorname{Tr}(\Hi)$. If $T$ is a non-negative, self-adjoint operator (i.e., $T = T^* $), then $\| T \|_{\operatorname{Tr}(\Hi)} = \operatorname{Tr}(T) = \sum_{i=1}^{\infty} \langle T u_i, u_i \rangle_{\operatorname{Tr}(\Hi)}$. Examples of non-negative, self-adjoint operators are covariance operators.
The three norms \eqref{eq:def-op-norm}--\eqref{def:trace_norm} satisfy
\begin{equation*} \label{eq:bounde_HS_trace_norm}
   \| \cdot \|_{\operatorname{op}} \le \| \cdot \|_{\operatorname{HS}(\Hi)} \leq \| \cdot \|_{\operatorname{Tr}(\Hi)}.
\end{equation*}

Let $(\mathbb{Y},\mathcal{A}, \mu)$ be a separable, $\sigma$-finite measure space. Then, the following isomorphisms hold, 
\begin{equation*}
L^{2}(\mathbb{Y} \times \mathbb{Y}, \mu \otimes \mu) \cong L^{2}(\mathbb{Y}, \mu) \otimes L^{2}(\mathbb{Y}, \mu) \cong L^{2}(\YY , \mu :L^{2}(\YY , \mu )).
\end{equation*}
Let $T \in L(\Hi)$ be a self-adjoint operator, and recall that a unitary operator $U$ is such that $U U^* = U^* U = I$, where $I$ is the identity operator. The spectral theorem for self-adjoint operators states that each self-adjoint operator is decomposable into a unitary operator and a multiplication operator; see Theorem 9.4.6 in \citet{comway1994course}. 
More precisely, there exist a measure space $(\YY,\mathcal{A}, \mu)$, a unitary operator $U$ and a multiplication operator $D_d$ on $(\YY,\mathcal{A}, \mu)$ associated to a bounded function $d$, such that 
\begin{equation} \label{eq:UNU=D}
UTU^{*}=D_d, \quad U:\Hi \to L^2(\mathbb{Y}), \quad D_d:L^2(\mathbb{Y}) \to L^2(\mathbb{Y}).
\end{equation}
Moreover, since $\Hi$ is separable, the measure space $(\YY,\mathcal{A}, \mu)$ is $\sigma$-finite; see Proposition 9.4.7 in \citet{comway1994course}. The multiplication operator $D_d$ is given by 
\begin{equation} \label{eq:D-def}
    D_d[f](s)\doteq d(s)f(s), \quad f \in  L^2(\mathbb{Y}), \; s \in \YY.
\end{equation}

\subsection{Properties of the linear process} \label{subsec:prelim-linear}

We collect here some general properties of linear processes with representation \eqref{eq:linear_process} as well as implications of imposing \eqref{eq:linear_process2}, where $T \in L(\Hi)$ is a self-adjoint operator. With regard to stochastic processes with sample paths in function spaces, we will make no distinction between the class of square integrable functions $\mathcal{L}^2$ and its respective family of equivalence classes $L^2$; see the first paragraph of Section 3 in \cite{Chara_Rack_Central}.

The linear process $\{X_n\}_{n \in \NN}$ 
in \eqref{eq:linear_process} converges $\P$-almost surely and in $L^2(\Omega:\Hi)$ if $\sum_{j=0}^\infty \| u_j \|^2_{\op} < \infty$, $\EE \varepsilon_0 = 0$ and $\EE \| \varepsilon_0 \|^2_{\Hi} < \infty$; see Lemma 7.1 in \cite{Bosq2000:Linear}. 
This implies that series that are SRD in the sense that $\sum_{j=0}^\infty \| u_j \|_{\op} < \infty$ converge.
Moreover, it shows that LRD series in the sense that $\sum_{j=0}^\infty \| u_j \|^{4/3}_{\op} < \infty$ convergence. 
However, for series that are LRD (in the sense of \eqref{eq:op-lrd}), both cases $\sum_{j=0}^\infty \| u_j \|^2_{\op} < \infty$ and $\sum_{j=0}^\infty \| u_j \|^2_{\op} = \infty$ are possible. For an instance of the latter case, note that with $\{ u_j \}$ satisfying \eqref{eq:uj-lrd}, we have $ \esssup_{s \in \YY} d(s) = 1/2$ if and only if $\sum_{j=0}^\infty \| u_j \|^2_{\op} = \infty$. In this case, the well posedness of the process is not given by the results in \cite{Bosq2000:Linear} and requires some additional assumptions. We first address the well-posedness of the linear processes $\{X_n\}_{n \in \ZZ}$ 
in \eqref{eq:linear_process} when $\sum_{j=0}^\infty \| u_j \|^2_{\op} = \infty$ for the case $\Hi = L^2(\YY)$. We then leverage the well-posedness on $L^2(\YY)$ and the spectral theorem \eqref{eq:UNU=D}, to derive the convergence of an $\Hi$-valued linear process \eqref{eq:linear_process}--\eqref{eq:linear_process2} for a general $\Hi$. 

First consider $\Hi = L^2(\YY)$ and rewrite the model \eqref{eq:linear_process} with \eqref{eq:uj-lrd} as
\begin{equation} \label{eq:XninL2}
X_{n}(r)=\sum_{j=0}^{\infty} (j+1)^{d(r)-1} \varepsilon_{n-j}(r), \quad r \in \YY,
\end{equation}
where $d(r) \in \left(0,\frac{1}{2}\right)$ for all $r \in \YY$ and $\{\varepsilon_{j}\}_{j \in \ZZ}$ is an $L^2(\mathbb{Y})$-valued i.i.d.\ sequence with
\begin{equation} \label{eq:sigmawithoutU}
\sigma(r,s)
\doteq	
\EE(\varepsilon_{0}(r) \varepsilon_{0}(s)),
\hspace*{0.2cm}
\sigma^2(r)
\doteq
\EE|\varepsilon_{0}(r)|^2,
\hspace*{0.2cm}
r,s \in \YY.
\end{equation}
The stochastic process $\{X_n\}$ in \eqref{eq:XninL2} (equivalently, \eqref{eq:linear_process} with \eqref{eq:uj-lrd}) has sample paths that belong to $L^2(\YY)$ a.s. if the following condition is satisfied
\begin{equation} \label{eq:integ-cond-well-def}
\int_{\YY} \frac{\sigma^2(r)}{1-2d(r)} \mu(dr) < \infty, \quad \text{implying that }
   \EE \| \varepsilon_0 \|^2_{L^2(\YY)} = \int_{\YY} \sigma^2(r) \mu(dr) < \infty ;
\end{equation}
see Proposition 3 in \cite{Chara_Rack_Central}. Note that \cite{Chara_Rack_Central} use a different notation $u_j = (j+1)^{-\tilde{d}(r)}$ with $\tilde{d}(r) = 1-d(r)$. In \eqref{eq:integ-cond-well-def}, their conditions are adjusted to our setting. 
Furthermore, the $L^2(\YY)$-valued series \eqref{eq:XninL2} converges in mean square and $\P$-almost surely if $d(s) < \frac{1}{2}$ $\mu$-a.s. under the second condition in \eqref{eq:integ-cond-well-def}; see Proposition 4 in \cite{Chara_Rack_Operator}, again after adjusting to the different notation. Hence Condition \eqref{eq:integ-cond-well-def} ensures the well-posedness of \eqref{eq:XninL2} when $\esssup_{s \in \YY} d(s) = 1/2$.

We now treat the case of a general $\Hi$. Invoking the spectral theorem, one can infer that \eqref{eq:linear_process} with \eqref{eq:linear_process2} also converges $\P$-a.s. if $d(s) < \frac{1}{2}$ $\mu$-a.s., where $d$ is the function associated to the multiplication operator $D_d$ in the decomposition of the spectral theorem \eqref{eq:UNU=D}--\eqref{eq:D-def}. We refer to p.\ 1445 in \cite{DUEKER} and calculations done in \eqref{eq:process_normal_perator_after_interchange}--\eqref{equality_process_Z} below.

Moreover, after accounting for the aforementioned change in notation, a CLT for the sample mean of $\{X_n\}_{n \in \ZZ}$ with values in $L^2(\YY)$ was shown in Proposition 4 of \citet{Chara_Rack_Central} under the conditions
    \begin{equation} \label{eq:cond-clt-mean}
        \int_{\YY} \frac{\sigma^2(r)}{d^2(r)} \mu(dr) < \infty, \quad \int_{\YY} \frac{\sigma^2(r)}{d(r)(1-2d(r))} \mu(dr) < \infty.
    \end{equation}
In the next section we compare these conditions with the corresponding ones for the weak convergence of the autocovariance operators.

For the sample autocovariances of $\{X_{n}\}_{n \in \ZZ}$ in \eqref{eq:XninL2} and its corresponding population quantities, we recast \eqref{eq:SampleAutoShort} as
\begin{equation} \label{eq:gammars}
\widehat{\gamma}_{N,h}(r,s) \doteq \frac{1}{N} \sum_{n=1}^{N} X_{n+h}(r)X_{n}(s)
\hspace{0.2cm}
\text{ and }
\hspace{0.2cm}
\gamma_{h}(r,s) \doteq \EE(X_{h}(r) X_{0}(s)).
\end{equation}
Note that henceforth, we use the notations $\widehat \Gamma_{N,h}, \Gamma_h$ for the autocovariance operators in a general space $\Hi$, and their lowercase counterparts $\widehat \gamma_{N,h}, \gamma_h$ for the special case $\Hi = L^2(\YY)$.

Finally, recall that, since our linear series admit second moments, the covariance operators are nuclear operators and therefore Hilbert-Schmidt; see p.\ 6 in \citet{Bosq2000:Linear}. 
This remains true for the empirical covariance operator that belongs $\P$-almost surely to the space $\Hi \otimes \Hi \cong \operatorname{HS}(\Hi)$; see pp. 36--37 in \cite{Bosq2000:Linear}.

\begin{remark} \label{re:slowly-varying-1}
\begin{enumerate}
    \item To ensure that \eqref{eq:linear_process2} indeed implies long-range dependence in the sense of \eqref{eq:op-lrd}, we employ the spectral decomposition \eqref{eq:UNU=D}. Then, \eqref{eq:op-lrd} holds if and only if $\esssup_{s \in Y} d(s) \ge 0$; see p.\ 1445 in \cite{DUEKER}.
    
    \item   Recall that, for $\Hi = \RR$, LRD is often modeled through a linear process with $u_j = j^{d-1}\ell(j)$ for $d \in (0,\frac{1}{2})$, where $\ell(j)$ is a slowly varying function; see Condition I on p.\ 17 in \cite{PipirasTaqqu}. The slowly varying function $\ell$ induces flexibility on the sequence $u_j$. While \eqref{eq:uj-lrd} naturally generalizes the real-valued model to the Hilbert space-valued setting allowing for a space-varying memory parameter, it is quite restrictive as a function in $j$. We emphasize that our results (Theorems \ref{thm1:srd} and \ref{th:main_result}) can be generalized to using $u_j = (j+1)^{D_d-I}\ell(j)$ with $\ell$ being a slowly varying function instead. Since the function $\ell$ is real-valued, one can adjust our proofs by incorporating the arguments for real-valued linear processes; see Chapter 2 in \cite{PipirasTaqqu}. For the sake of 
  clarity, we only consider here the case $\ell(j) = 1$ and focus on the remaining technical difficulties.
\end{enumerate}
\end{remark}

\section{Main Results}
\label{s:mainresults}

In this section, we introduce the limiting objects and state the convergence results for the sample autocovariance operators \eqref{eq:SampleAutoShort}. 
From here on we distinguish the following two cases, roughly corresponding to the two regimes $d(s) \in \left(0,\frac{1}{4}\right)$ and $d(s) \in \left(\frac{1}{4},\frac{1}{2}\right)$ for each $s \in \mathbb{Y}$.
First, suppose that $\{u_j\}_{j \in \NN_0}$ is such that
\begin{equation} \label{cond:srd-improved_1}
\sum_{j=0}^\infty \| u_j \|_{\operatorname{op}}^{4/3} < \infty,
\end{equation}
which in the following we refer to as \textit{first regime}. Second, suppose that $T$ is a self-adjoint operator and that
\begin{equation} \label{eq:second_regime}
    u_{j} = (j + 1)^{T-I},
    \hspace{0.2cm}
    T = UD_dU^* 
    \hspace{0.2cm}
    \text{ with }
    \hspace{0.2cm}
     d(s) \in \left(\frac{1}{4},\frac{1}{2}\right),   \; s \in \YY,
\end{equation}
which we refer to from now on as \textit{second regime}. In \eqref{eq:second_regime}, $U,D_d$ are as in the spectral theorem \eqref{eq:UNU=D}. In particular, under the second regime, the series in \eqref{cond:srd-improved_1} diverges.

We start by defining the covariance operator of the limiting Gaussian process in the first regime. We view the limiting Gaussian element for a single time lag as an element of $\operatorname{HS}(\Hi)$. For the joint convergence of autocovariance operators across time lags $0,1,\dots,H$, for $H \in \NN$, the corresponding Gaussian limit is an element of $\operatorname{HS}(\Hi)^{\times (H+1)}$, and its covariance operator is an element of trace class on the space $\operatorname{HS}(\Hi)^{\times (H+1)}$. It can be identified with an operator taking the block form
\begin{equation} \label{eq:Sigma}
    \Sigma \doteq (\Sigma^{(p,q)}_\Gamma)_{p,q = 0,\dots,H}.
\end{equation}
For $p,q = 0,\dots,H$, and $T \in \operatorname{HS}(\Hi)$, the cross-covariance operator $\Sigma^{(p,q)}_\Gamma$ is given by 
\begin{equation} \label{eq:covar-oper-Sigma}
    \Sigma_\Gamma^{(p,q)}(T) = \sum_{m = 0}^\infty  \Gamma_{m+p-q} T  \Gamma_{m} + \sum_{m = 0}^\infty  \Gamma_{m+q} T \Gamma_{m-p} +  A_q(\Lambda - \Phi) A_p (T),
\end{equation}
where, upon recalling that $\langle \varepsilon_0, \cdot \rangle_{\Hi} \varepsilon_0 \in \operatorname{HS}(\Hi)$, we have $\Lambda, \Phi \in \operatorname{HS} \left(\operatorname{HS}(\Hi) \right)$ with
\begin{equation} \label{def:def:Lambda-Phi}
    \Lambda(T) \doteq \EE \left( \langle  \langle \varepsilon_0, \cdot \rangle_{\Hi}  \varepsilon_0, T \rangle_{\operatorname{HS(\Hi)}} ( \langle \varepsilon_0, \cdot \rangle_{\Hi} \varepsilon_0)   \right), \quad \Phi(T) \doteq \langle C, T + T^* \rangle_{\operatorname{HS(\Hi)}} C + \langle C , T \rangle_{\operatorname{HS(\Hi)}} C, 
\end{equation}
and $C, \Gamma_h \in \operatorname{HS}(\Hi) , A_p \in \operatorname{HS}\left(\operatorname{HS}(\Hi) \right)$ with
\begin{equation} \label{def:def:Aq-C}
  C \doteq \EE \left( \langle \varepsilon_0, \cdot \rangle_{\Hi} \varepsilon_0 \right), \quad \Gamma_h \doteq \EE \langle X_h, \cdot \rangle_{\Hi} X_0, \quad A_p(T) \doteq \sum_{j \in \NN} u_{j+p} T u_j^*.
\end{equation}
Note that the two representations of  $\Gamma_h$ in \eqref{eq:SampleAutoShort} and \eqref{def:def:Aq-C} are due to the aforementioned isomorphism $\Hi \otimes \Hi \cong \operatorname{HS}(\Hi)$. 
The representation \eqref{eq:covar-oper-Sigma}--\eqref{def:def:Aq-C} is due to Lemma 3 in \cite{Mas}, but the notation is slightly changed to be consistent with the present paper.

\begin{theorem}[First regime] \label{thm1:srd}
    Let $\{X_n\}_{n \in \ZZ}$ be an $\mathbb{H}$-valued linear process \eqref{eq:linear_process}, and consider its autocovariance operators $(\widehat \Gamma_{N,h}, h = 0,\dots,H)$ given in \eqref{eq:SampleAutoShort} with $\widehat \Gamma_{N,h} \in \operatorname{HS}(\Hi)$ for all $h$ with a slight abuse of notation. Suppose \eqref{cond:srd-improved_1} and 
\begin{equation} \label{eq:Th1-lemma-ass1}
  \EE \| \varepsilon_{0} \|^4_{\mathbb{H}}< \infty. 
\end{equation}
    Then,
    \begin{equation}
        \sqrt{N} \begin{pmatrix}
            \widehat \Gamma_{N,0} - \Gamma_{0} \\
            \vdots \\
            \widehat \Gamma_{N,H} - \Gamma_{H}
        \end{pmatrix} \xrightarrow{d} 
        G \doteq
        \begin{pmatrix}
            G_0 \\
            \vdots \\
            G_H
        \end{pmatrix} \in (\operatorname{HS}(\mathbb{H}))^{\times (H+1)},
    \end{equation}
    where the weak convergence holds in the topology of $(\operatorname{HS}(\mathbb{H}))^{\times (H+1)}$, and $G$ is the centered Gaussian element of $\operatorname{HS}(\Hi)^{\otimes (H+1)}$ with covariance operator $\Sigma$ given in \eqref{eq:Sigma}--\eqref{def:def:Aq-C}.
\end{theorem}

The proof of Theorem \ref{thm1:srd} is postponed to Section \ref{subs:lemma0}. Note that Theorem \ref{thm1:srd} does not require any structural assumptions on the sequence $\{u_j\}_{j\in \NN_0}$ besides \eqref{cond:srd-improved_1}. In particular, \eqref{eq:linear_process2} is not needed. The following corollary states the implications of Theorem \ref{thm1:srd} on the linear series with values in $L^2(\mathbb{Y})$ with \eqref{eq:linear_process2}, as written in \eqref{eq:XninL2}.

\begin{corollary}\label{cor:srd-L^2}
    Let $\{X_n\}_{n \in \ZZ}$ be an $L^2(\mathbb{Y})$-valued process \eqref{eq:XninL2} and consider its autocovariance operators $( \widehat \gamma_{N,h}, h = 0,\dots,H)$ given in \eqref{eq:gammars}, where we consider $\widehat \gamma_{N,h}$ as an element of $\operatorname{HS}(L^2(\mathbb{Y}))$ for all $h$, with a slight abuse of notation. Suppose $\esssup_{s \in \YY} d(s) < \frac{1}{4}$ and
\begin{equation} \label{eq:Th2-lemma-ass1}
  \EE \| \varepsilon_{0} \|^4_{L^2(\mathbb{Y})} < \infty. 
\end{equation}
Then, 
    \begin{equation}
        \sqrt{N} \begin{pmatrix}
            \widehat \gamma_{N,0} - \gamma_{0} \\
            \vdots \\
            \widehat \gamma_{N,H} - \gamma_{H}
        \end{pmatrix} \xrightarrow{d} 
        G \doteq
        \begin{pmatrix}
            G_0 \\
            \vdots \\
            G_H
        \end{pmatrix} \in \left(\operatorname{HS}(L^2(\mathbb{Y})\right)^{\times (H+1)},
    \end{equation}
    where the convergence holds in the topology of $(\operatorname{HS}(L^2(\mathbb{Y})))^{\times (H+1)}$, $G$ is a centered Gaussian element,
     and the covariance operator of $G$ is the operator $\Sigma$ given in \eqref{eq:Sigma}--\eqref{def:def:Aq-C}, for $p,q=0,\dots, H$, and with $\Hi$ replaced by $L^2(\YY)$.
\end{corollary}

\begin{proof}
 Let $\esssup_{s \in \YY} d(s) = \frac{1}{4} - \delta$, for some $\delta \in (0,\frac{1}{4})$. We see that
    \begin{equation*}
        \sum_{j=0}^{\infty} \|u_j\|^{4/3}_{\operatorname{op}} = 
        \sum_{j=0}^{\infty} \left((j+1)^{\esssup_{s \in \mathbb{Y}} d(s) - 1} \right)^{4/3} = 
        \sum_{j=0}^{\infty} (j+1)^{-1 - \frac{4 \delta}{3}} < \infty.
    \end{equation*}
    Then, Corollary \ref{cor:srd-L^2} follows from Theorem \ref{thm1:srd}.
\end{proof}

The second regime is more challenging. We first consider fluctuations for the autocovariance operators defined in \eqref{eq:gammars} for $X_n$ taking values in $\Hi = L^2(\YY,\mu)$ (Theorem \ref{le:limit_theorem_L2}). Then, we leverage this result and the spectral theorem for self-adjoint operators to identify the fluctuations of the autocovariance operator defined in \eqref{eq:SampleAutoShort} (Theorem \ref{th:main_result}). 

In the second regime, the limit is no longer Gaussian. Instead, the resulting limit process can be represented as a double Wiener-It\^o integral with sample paths in $L^2(\mathbb{Y}^2)$, applied to a certain kernel $f \in L^2(\YY^2)$ defined in \eqref{equality_f^(u,l)_t,(m,s)} below. This limiting object is a generalization of the Rosenblatt distribution, allowing for a continuum of memory parameters, and reflecting the fact that the underlying process is infinite-dimensional. Double Wiener-It\^o integrals with values in function spaces are introduced and defined in Section \ref{se:biW} below.

We define the suitable scaling as a multiplication operator on $L^2(\mathbb{Y}^2:\RR^{H+1})$ by
\begin{equation} \label{eq:Xi-def}
    \Xi_N[f](r,s) \doteq N^{1-d(r) - d(s)} f(r,s), \quad f \in L^2(\mathbb{Y}^2:\RR^{H+1}), \quad r , s \in \YY.
\end{equation}

\begin{theorem} \label{le:limit_theorem_L2}
Let $\{X_{n}\}_{n \in \ZZ}$ be an $L^2(\mathbb{Y})$-valued linear process \eqref{eq:XninL2} and consider its autocovariance operators $(\widehat \gamma_{N,h}, h = 0,\dots,H)$ given in \eqref{eq:gammars}. Suppose $ d(s) \in \left(\frac{1}{4},\frac{1}{2}\right)$, 
\begin{equation} \label{eq:lemma-ass1}
  \EE \| \varepsilon_{0} \|^4_{L^4(\mathbb{Y})} < \infty. 
\end{equation}
and
\begin{equation} \label{eq:lemmaint2}
\int_{\mathbb{Y}} \frac{\sigma^2(r)}{(1-2d(r))(2d(r)-1/2)}	\mu(dr)  < \infty.
\end{equation}
Then, 
\begin{equation*}
\Xi_N \begin{pmatrix}
    \widehat{\gamma}_{N,0}-\gamma_{0} \\
    \vdots 
    \\
    \widehat{\gamma}_{N,H}-\gamma_{H}
\end{pmatrix} 
\overset{\distr}{\to}
\begin{pmatrix}
    \mathfrak{R} \\
    \vdots \\
     \mathfrak{R}
\end{pmatrix}
 \in L^2(\mathbb{Y}^2:\RR^{(H+1)}),
\end{equation*}
where $\mathfrak{R}$ is given in \eqref{equality_f^(u,l)_t,(m,s)} below and $\Xi_N$ in \eqref{eq:Xi-def}.
\end{theorem}

We continue with some remarks on the above assumptions.

\begin{remark}
\begin{enumerate}
    \item Since $(\YY,\mu)$ is not necessarily a finite measure space, the condition $\EE \|\varepsilon_0\|^4_{L^4(\YY)} < \infty$ does not imply $\EE \|\varepsilon_0\|^2_{L^2(\YY)} < \infty$. However $\EE \|\varepsilon_0\|^2_{L^2(\YY)} < \infty$ is contained within Condition \eqref{eq:lemmaint2}. 

    \item If there exists some $\delta > 0$ such that $1/4 + \delta < d(r) < \frac{1}{2} - \delta$ for all $r \in \YY$, then the Conditions $\EE \|\varepsilon_0\|^2_{L^2(\YY)} < \infty$ and \eqref{eq:lemmaint2} are equivalent.

    \item The estimates in \eqref{eq:lemma-ass1}--\eqref{eq:lemmaint2} ensure that \eqref{eq:integ-cond-well-def} is satisfied, hence the process $\{X_n\}_{n \in \ZZ}$ has a.s. sample paths in $L^2(\YY)$; see Section \ref{subsec:prelim-linear}. Moreover, the conditions in \eqref{eq:cond-clt-mean} are also satisfied due to \eqref{eq:lemma-ass1}--\eqref{eq:lemmaint2} and since $d(r) > \frac{1}{4}$. Therefore, following upon our discussion in Section \ref{subsec:prelim-linear}, the conditions of Theorem \ref{le:limit_theorem_L2} imply that the sample mean of $\{ X_n\}_{n \in \ZZ}$ satisfies a CLT.
    
\end{enumerate}
\end{remark}

We turn to our final main result, which lifts the second regime to linear processes taking values in general Hilbert spaces. We start by introducing the right scaling. Define, for some unitary operator $U: \Hi \to L^2(\YY)$, the operator $\Delta^U_N \in L\left( (\Hi \otimes \Hi)^{\times (H+1)} \right)$ by
\begin{equation} \label{eq:scaling:Delta}
\Delta^U_{N}[f] \doteq 
 \left( U^* \otimes U^* \right)^{\times (H+1)} \Xi_N \left( U \otimes U \right)^{\times (H+1)} [f] , \quad f \in (H \otimes H)^{\times(H+1)},
\end{equation}
where $\Xi_N$ is the operator defined in \eqref{eq:Xi-def}. We then have the following theorem.

\begin{theorem}[Second regime] \label{th:main_result}
Let $\{X_{n}\}_{n \in \ZZ}$ be an $\Hi$-valued linear process \eqref{eq:linear_process}, \eqref{eq:linear_process2}, and consider its autocovariance operators $(\widehat \Gamma_{N,h}, h = 0,\dots,H)$ given in \eqref{eq:SampleAutoShort}. Suppose \eqref{eq:second_regime},
\begin{equation} \label{eq:thm-lrd-ass1}
  \EE \|U \varepsilon_{0} \|^4_{L^4(\mathbb{Y})} < \infty 
\end{equation}
and
\begin{equation} \label{eq:thm-lrd-ass2}
\int_{\mathbb{Y}} \frac{\EE (U \varepsilon_0)^2(r)}{(1-2d(r))(2d(r) - 1/2)}	\mu(dr)  < \infty.
\end{equation}
Then,
\begin{equation*}
\Delta^U_{N} \begin{pmatrix}
    \widehat{\Gamma}_{N,0}-\Gamma_{0} \\
    \vdots \\
    \widehat{\Gamma}_{N,H}-\Gamma_{H}
\end{pmatrix}
\overset{\distr}{\to}
\begin{pmatrix}
\mathfrak{Z}_U \\
\vdots \\
\mathfrak{Z}_U
\end{pmatrix}
\in (\Hi \otimes \Hi)^{\times (H+1)},
\end{equation*}		
where $\Delta^U_N$ is defined in \eqref{eq:scaling:Delta} and $\mathfrak{Z}_U$ is defined in \eqref{equality_f^(u,l)_t,(m,s)}--\eqref{eq:Ztilde} below.
\end{theorem}

\section{Double Wiener-It\^o Integrals in Function Spaces}
\label{se:biW}

In this section, we construct double Wiener-It\^o integrals with values in $L^2(\mathbb{Y}^2)$, extending the work of \cite{Nor94} in the direction of integrators with spatial dependence, and the work of \cite{FoxTaqqu87} in the direction of double integrals with sample paths in an infinite dimensional Hilbert space. Our candidate Wiener-It\^o integrals are defined by integrating with respect to a family of dependent measures $\{W^{(r)}\}_{r \in \YY}$. In order to define a double Wiener-It\^o  integral in $L^2(\mathbb{Y}^2)$, we first define the integral for special kernels, and then use an approximation of $f$ in terms of such special kernels. Note that $L^2(\mathbb{Y}^2)$ is a complete, separable, $\sigma$-finite measure space. Our construction is applicable for multiple Wiener-It\^o integrals of any order, but we restrict ourselves to double Wiener-It\^o integrals for simplicity. We start by recalling some elementary notions of the usual, $\RR$-valued double Wiener-It\^o integrals.

We write $(\RR,\mathcal{B},\lambda)$ for the measure space with the usual Borel topology $\mathcal{B}$ and the Lebesgue measure $\lambda$. For every $\psi \in L^2(\RR^2, \lambda^2)$, the (standard) double Wiener-It\^o integral with respect to a Gaussian random measure $G$ is defined by
\[
I_2(\psi) \doteq \int_{\RR^2} \psi(x_1, x_2) \, G(dx_1) G(dx_2). 
\]
A Wiener-It\^o integral of $\psi \in L^2(\RR^2, \lambda^2)$ with regard to dependent integrators \eqref{eq:tilde-G-r-s} was constructed in \cite{FoxTaqqu87}.
Let $G_1,G_2$ be two dependent Gaussian random measures satisfying
\begin{equation} \label{eq:tilde-G-r-s}
    \EE G_1 (A) G_2(B) = \sigma_{G_1 G_2} \lambda(A \cap B),
    \hspace{0.2cm}
    A,B \in \mathcal{B},
\end{equation}
for some $\sigma_{G_1 G_2}>0$. Then, \cite{FoxTaqqu87} introduced the corresponding double Wiener-It\^o integral
\begin{equation} \label{eq:tilde-I-2-def-G}
    \tilde I_2(\psi) = \int_{\RR^2} \psi(x_1, x_2) \, G_1(dx_1) G_2(dx_2), \quad \psi \in L^2(\RR^2). 
\end{equation}

In contrast to \cite{FoxTaqqu87}, we fix a correlation function $\sigma: \YY^2 \to \RR$, denote $\sigma^2(r) \doteq \sigma(r,r)$, and consider a family of dependent Gaussian random measures $\{W^{(r)}(A): A \in \mathcal{B}_0 \}_{r \in \YY}$, where $\mathcal{B}_0 = \{A \in \mathcal{B}: \nu^{(r)}(A) = \sigma^2(r) \lambda(A) < \infty, \text{for all } r \in \YY \}$. All measures are defined on a common probability space $(\Omega, \mathcal{F}, \P)$. Let their covariances be given by, for $r,s \in \YY$, and $x \neq y$,
\begin{equation} \label{eq:covWW}
\EE W^{(r)}(dx)=0,
\hspace{0.2cm}
\EE W^{(r)}(dx) W^{(s)}(dx)=\sigma(r,s)dx,	\hspace*{0.2cm}
\EE W^{(r)}(dx) W^{(s)}(dy)=0.
\end{equation}
Note that, for fixed $r$, the measure $W^{(r)}(dx)$ is viewed as the increment of an infinitesimal interval $dx$ of a Brownian
motion with coefficient $\sigma(r) \doteq \sqrt{\sigma^2(r)}$. For fixed $r,s \in \YY$, we set
\begin{equation} \label{eq:tilde-I-2-def}
    \tilde I^{(r,s)}_2(\psi) = \int_{\RR^2} \psi(x_1, x_2) \, W^{(r)}(dx_1) W^{(s)}(dx_2), \quad \psi \in L^2(\RR^2),
\end{equation}
where $\tilde I_2$ is the stochastic integral defined in \eqref{eq:tilde-G-r-s}--\eqref{eq:tilde-I-2-def-G}, with $G_1 = W^{(r)}, G_2 = W^{(s)}$. 

We aim to integrate with respect to the family of measures $\{W^{(r)}(A): A \in \mathcal{B}_0 \}_{r \in \YY}$. Therefore, we define the set of admissible kernels for such integrals. We call a map $f^{(r,s)}(x,y) \doteq f(r,s,x,y), \; (r,s,x,y) \in \YY^2 \times \RR^2$ an \textit{admissible kernel} if (i) $f^{(r,s)}(x,y)$ is a $
\mu^{\otimes 2} \otimes \lambda^{\otimes 2}$-measurable function, (ii) the family $\{\tilde I^{(r,s)}_2(f^{(r,s)}): r,s \in \YY\}$ viewed as a stochastic process is measurable, where for fixed $r,s \in \YY$, $\tilde I^{(r,s)}_2(f^{(r,s)})$ is defined as in \eqref{eq:covWW}--\eqref{eq:tilde-I-2-def}, and (iii) the map from $\YY^2 \times \RR^2$ to $\RR$ defined by 
\begin{equation} \label{eq:defn-f-hat}
    (r,s,x,y) \mapsto f^{\sigma}(r,s,x,y) \doteq f^{(r,s)}(x,y) \sigma (r) \sigma(s)
\end{equation}
is in $L^2(\RR^2: L^2(\YY^2))$ (with a slight abuse of notation).
We denote by $\mathcal{H}^2 \doteq \mathcal{H}^2(\YY^2 \times \RR^2)$ the space of admissible kernels satisfying Conditions (i)--(iii). Note that the space $\mathcal{H}^2$ depends on the selection of the function $\sigma$, but we suppress this dependence for notational simplicity. We present some remarks about the space $\mathcal{H}^2$.

\begin{remark} \label{rmk:meas}
\begin{enumerate}
    \item If $\YY$ is a separable metric space equipped with its Borel $\sigma$-algebra and if the kernel $f$ is a jointly measurable function on $\YY^2 \times \RR^2$, then $\{\tilde I^{(r,s)}_2 : r, s \in \YY\}$ has a measurable modification and can hence be assumed to be measurable; see \cite{Nor94}, p.\ 338. 
    \item In particular, Condition (iii) implies that $f^{(r,s)}(\cdot,\cdot) \in L^2(\RR^2, \lambda^{\otimes 2})$, for $\mu$-a.e. $r,s \in \YY$ since $\sigma^2(r) > 0$ for all $r \in \YY$.
\end{enumerate}
\end{remark}


We proceed with defining our double Wiener-It\^o integral for \textit{special kernels}. We call an admissible kernel $f \in \mathcal{H}^2$ \text{special}, if there exists an $N \in \NN$ and a system $\{\Delta_1, \dots, \Delta_N\}$ of disjoint sets in $\mathcal{B}_0$ such that, for $x,y \in \RR$, $r,s \in \YY$,
\begin{equation} \label{eq:simple-kernel}
f^{(r,s)}(x,y) = \begin{cases}
    c_{i_1, i_2}(r,s)  \mathds{1}_{\Delta_{i_1}}(x)  \mathds{1}_{\Delta_{i_2}}(y) & \text{for some }i_1,i_2 \in \{1,\dots,N\} \text{ such that } i_1 \neq i_2,   \\
    0 & \text{else},
\end{cases} 
\end{equation}
where $c_{i_1, i_2} \in L^2(\YY^2,\mu^{\otimes 2})$ for all $(i_1, i_2) \in \{1, \dots, N\}^{\times 2}$. Then, the special admissible kernel $f^{(\cdot,\cdot)}(x,y)$ takes values in $L^2(\YY^2,\mu^{\otimes 2})$, and we can define the $L^2(\YY^2,\mu^{\otimes 2})$-valued stochastic process $\mathcal{I}_2(f)$ by
\[
\mathcal{I}_2(f)(r,s) \doteq \sum_{\substack{i_1,  i_2 = 1 \\ i_1 \neq i_2}}^N c_{i_1, i_2}(r,s) W^{(r)}(\Delta_{i_1}) W^{(s)}(\Delta_{i_2}), \quad r,s \in \YY.
\]

\begin{definition} \label{def:double-in-banach}
Let $(\YY,\mu)$ be a $\sigma$-finite measure space, and let $\{W^{(r)}\}_{r \in \YY}$ be a family of Gaussian random measures on $(\RR, \lambda)$ with covariances \eqref{eq:covWW}. Consider a kernel $f: \YY^2 \times \RR^2 \to \RR$. We say that there exists an $L^2(\YY^2)$-valued double Wiener-It\^o integral of the kernel $f$, denoted $\mathcal{I}_2(f)$, if there exists a sequence $\{f_{(n)}: n \geq 1\}$ of $L^2(\YY^2)$-valued special kernels, defined on $\YY^2 \times \RR^2$, such that
\begin{enumerate}[label=(\roman*)]
    \item $\lim_{n \to \infty} \|f^{(r,s)} - f_{(n)}^{(r,s)}\|_{L^2(\RR^2)} = 0$ $\mu^{\otimes 2}$-a.s.,
    \label{item:Def42_1}
    \item the sequence $\{\mathcal{I}_2(f_{(n)}): n \geq 1\}$ is a Cauchy sequence in $L^1(\Omega:L^2(\YY^{2}))$.
    \label{item:Def42_2}
\end{enumerate}
In that case, we define the double Wiener-It\^o integral by 
\[
\mathcal{I}_2(f) \doteq \lim_{n \to \infty} \mathcal{I}_2(f_{(n)}),
\]
where the limit is again understood in $L^1(\Omega:L^2(\YY^{2}))$. For an $\RR^d-$valued kernel $f = (f_1,\dots,f_d)$ such that $\mathcal{I}_2(f_i)$ exists for all $i = 1,\dots,d$, we fix the notation $\mathcal{I}_2(f) \doteq (\mathcal{I}_2(f_1),\dots,\mathcal{I}_2(f_d))$.
\end{definition} 

Whenever the underlying covariance structure of $\{W^{(r)}\}_{r \in \YY}$ is clear from the context, we write $\mathcal{I}_2$; otherwise, we write $\mathcal{I}_2^\sigma$ to emphasize that the dependence of the underlying Gaussian random measure is characterized through \eqref{eq:covWW}.

\begin{remark} \label{rmk:int-taqqu-ours}
It follows that if there exists an $L^2(\YY^2)$-valued double Wiener-It\^o integral $\mathcal{I}_2(f)$ of a kernel $f$, then the double Wiener-It\^o integral $\mathcal{I}_2(f)$ has $\P$-a.s. sample paths in $L^2(\YY^2)$ and is well defined, i.e., the definition of $\mathcal{I}_2(f)$ does not depend on the choice of a sequence $\{f_{(n)}: n \geq 1\}$ of $L^2(\YY^2)$-valued special functions. Moreover, the double Wiener-It\^o integral $\mathcal{I}_2$ coincides with the one defined in \cite{FoxTaqqu87}, i.e.,
\begin{equation}
    \tilde I^{(r,s)}_2(f^{(r,s)}) = \mathcal{I}_2(f)(r,s), \quad \mu^{\otimes 2}\text{-a.s.},
\end{equation} 
where $\tilde I_2$ is the stochastic integral defined in \eqref{eq:covWW}--\eqref{eq:tilde-I-2-def}.
 \end{remark}

We show that the $L^2(\YY^2)$-valued double Wiener-It\^o integral we defined exists for all kernels $f \in \mathcal{H}^2$. Note that the abstract machinery employed by \cite{Nor94} to define such multiple stochastic integrals with values in functional spaces is not readily available. The analysis in \cite{Nor94} relies heavily on a hypercontractivity property of the double stochastic integrals. However, we are not aware of whether the  hypercontractivity property remains valid when the integration is conducted with regard to Gaussian measures with spatial dependence.

\begin{theorem} \label{lemma:simple-dominated-above}
    Let $f \in \mathcal{H}^2$. Then, the $L^2(\YY^2)$-valued double Wiener-It\^o integral $\mathcal{I}_2(f)$ of the kernel $f$ exists and
    \begin{equation} 
         \EE \left( \int_{\YY^2} |\mathcal{I}_2(f^{(r,s)})|^2 \mu(dr) \mu(ds)\right)^{1/2} < \infty.
    \end{equation}
    Moreover, the sequence $\{f_{(n)}\}_{n \in \NN}$ approximating $f \in \mathcal{H}^2(\YY^2 \times \RR^2)$ can be selected such that, for all $n \in \NN$,
    \begin{equation} \label{eq:fn-dom}
        f_{(n)}^{(r,s)}(x,y) \le f^{(r,s)}(x,y), \quad \text{for all } (x,y) \in \RR^2, (r,s) \in \YY^2.
    \end{equation}
\end{theorem}

\begin{proof}
    We must show that there exists a sequence $f^{(r,s)}_{(n)}$ such that the definition requirements are satisfied for all $r,s \in \YY$. We start with item \ref{item:Def42_1} of Definition \ref{def:double-in-banach}.

    For $r,s\in \YY$ fixed, consider the map $(x,y) \mapsto f^{(r,s)}(x,y)$ as an element of $L^2(\RR^2)$. The set of special kernels is dense in $L^2(\RR)$ (for all $r,s$), and so there exists a family of functions $\{f_{(n)}^{(r,s)}\}_{n \in \NN, r,s \in \YY}$ such that the first item \ref{item:Def42_1} of Definition \ref{def:double-in-banach} is true.

    Moreover, we see that the simple functions can be chosen such that \eqref{eq:fn-dom} is true as follows. We define a new sequence $\{\tilde f_{(n)}\}_{n \in \NN}$ by
    \begin{equation*}
        \tilde f_{(n)}^{(r,s)}(x,y) \doteq
        \begin{cases}
            f_{(n)}^{(r,s)}(x,y) & \text{if } f_{(n)}^{(r,s)}(x,y) \le f^{(r,s)}(x,y), \\
            f^{(r,s)}(x,y) - (f_{(n)}^{(r,s)}(x,y) -f^{(r,s)}(x,y)) & \text{if } f_{(n)}^{(r,s)}(x,y) > f^{(r,s)}(x,y).
        \end{cases}
    \end{equation*}
    Then, $f_{(n)}^{(r,s)}(x,y) \le f^{(r,s)}(x,y)$ for all $x,y \in \RR, r,s \in \YY$. Moreover, for fixed $r,s \in \YY$,
    \begin{align*}
        \|f^{(r,s)} - \tilde f^{(r,s)}_{(n)} \|^2_{L^2(\RR^2)} 
        &= \int_{\RR^2}
        | f^{(r,s)}(x,y) - \tilde f^{(r,s)}_{(n)}(x,y)|^2 dx dy
        \\&=  \int_{\RR^2}
        | f^{(r,s)}(x,y) - f^{(r,s)}_{(n)}(x,y)|^2 dx dy \to 0,
    \end{align*}
as $n \to \infty$.

We now turn to item \ref{item:Def42_2} of Definition \ref{def:double-in-banach}. We show that $\{\mathcal{I}_2(f_{(n)})\}_{n \in \NN}$ is Cauchy in the complete metric space $L^1(\Omega:L^2(\YY^2))$. In view of Remark \ref{rmk:int-taqqu-ours}, we write for $m,n \in \NN$ and special kernels $f_{(m)},f_{(n)}$,
\begin{align}
    &
    \| \mathcal{I}_2(f_{(m)}) -  \mathcal{I}_2(f_{(n)})\|_{L^1(\Omega:L^2(\YY^2))} \nonumber
    \\&= \EE \left| \int_{\YY^2} \left| \tilde I^{(r,s)}_2(f_{(n)}^{(r,s)}) -  \tilde I^{(r,s)}_2(f_{(m)}^{(r,s)}) \right|^2 \mu(dr) \mu(ds) \right|^{1/2} \nonumber\\
    &\le \left| \EE \int_{\YY^2} \left| \tilde I^{(r,s)}_2(f_{(n)}^{(r,s)}) -  \tilde I^{(r,s)}_2(f_{(m)}^{(r,s)}) \right|^2 \mu(dr) \mu(ds) \right|^{1/2} \nonumber\\
    &\le \bigg( \int_{\YY^2} \int_{\RR^4} \left(f_{(m)}^{(r,s)}(x_1,x_2) - f_{(n)}^{(r,s)}(x_1,x_2)  \right) \left( f_{(m)}^{(r,s)}(y_1,y_2) - f_{(n)}^{(r,s)}(y_1,y_2) \right)  \nonumber\\ 
    &\quad \quad \times \EE (W^{(r)}(dx_1) W^{(s)}(dx_2) W^{(r)}(dy_1) W^{(s)}(dy_2) ) \mu(dr) \mu(ds)  \bigg)^{1/2},
    \label{al:pp_line4}
\end{align}
where the third line follows from Jensen's inequality by exchanging expectation and the square root operation. Recall that in this relation, with the integration over $\RR^4$ we are excluding the diagonals $x_1=x_2, y_1 = y_2$. From the calculation in \eqref{al:pp_line4}, we have that
\begin{align}
     &\| \mathcal{I}_2(f_{(m)}^{(r,s)}) -  
     \mathcal{I}_2(f_{(n)}^{(r,s)})\|^2_{L^1(\Omega:L^2(\YY^2))} \nonumber\\
     &\le   \int_{\YY^2} \int_{\RR^2} (f_{(m)}^{(r,s)}(x_1,x_2) - f_{(n)}^{(r,s)}(x_1,x_2))^2 \sigma^2(r) \sigma^2(s) dx_1 dx_2  \mu(dr) \mu(ds)  \nonumber\\
    &\quad \quad  + \int_{\YY^2} \int_{\RR^2} \left( f_{(m)}^{(r,s)}(x_1,x_2) - f_{(n)}^{(r,s)}(x_1,x_2)  \right) \left( f_{(m)}^{(r,s)}(x_2,x_1) - f_{(n)}^{(r,s)}(x_2,x_1) \right) 
    \nonumber
    \\
    &\quad \quad \times \sigma^2(r,s)  dx_1 dx_2  \mu(dr) \mu(ds) \nonumber \\
    &\le 4\int_{\YY^2} \int_{\RR^2} (f_{(m)}^{(r,s)}(x_1,x_2) - f_{(n)}^{(r,s)}(x_1,x_2))^2 \sigma^2(r) \sigma^2(s) dx_1 dx_2  \mu(dr) \mu(ds) ,
    \nonumber
\end{align}
where the last line follows from two iterations of Cauchy-Schwarz. Moreover, note that
\begin{multline}
     \int_{\YY^2} \int_{\RR^2} (f_{(m)}(x_1,x_2) - f_{(n)}(x_1,x_2))^2 \sigma^2(r) \sigma^2(s) dx_1 dx_2  \mu(dr) \mu(ds)  \\
    \le 4 \int_{\YY^2} \int_{\RR^2} (f^{(r,s)}(x_1,x_2))^2 \sigma^2(r) \sigma^2(s) dx_1 dx_2  \mu(dr) \mu(ds) 
    = 4 \|f^\sigma\|_{L^2(\RR^2: L^2(\YY^2))}< \infty,
\end{multline}
since $f \in \mathcal{H}^2$ and with $f^\sigma$ defined in \eqref{eq:defn-f-hat}. By DCT, it follows that we can select $n,m \in \NN$ large enough such that
\begin{equation}
    \| \mathcal{I}_2(f_{(m)}) -  \mathcal{I}_2(f_{(n)})\|^2_{L^1(\Omega:L^2(\YY^2))} \le \varepsilon.
\end{equation}
The same argument shows that, for general $f \in \mathcal{H}^2$ 
\begin{equation}
\EE \left( \int_{\YY^2} |\mathcal{I}_2(f^{(r,s)})|^2 \mu(dr) \mu(ds)\right)^{1/2} < \infty,
\end{equation}
which concludes the proof.
\end{proof}

\subsection{Rosenblatt distribution with sample paths in $L^2(\YY^2)$}
We define, for $d(r) \in \left(\frac{1}{4},\frac{1}{2}\right), r \in \YY$,
\begin{equation} \label{equality_f^(u,l)_t,(m,s)}
\mathfrak{f}^{(r,s)}(x_{1},x_{2}):=
\int_{0}^{1} (v-x_{1})_{+}^{d(r)-1} (v-x_{2})_{+}^{d(s)-1} dv, \quad \mathfrak{R} \doteq \mathcal{I}_2(\mathfrak{f}),
\end{equation}
where $x_{+} \doteq \max\{0,x\}$. We say that $\mathfrak{R}$ follows the \textit{Rosenblatt distribution with sample paths in} $L^2(\YY^2)$. This aligns with the terminology presented in \cite{veillette2013}, who consider the Rosenblatt distribution in the real-valued case.

To ensure that $\mathcal{I}_2(\mathfrak{f})$ is well-defined, we require the following result.

\begin{lemma}
    Let $d(r) \in \left(\frac{1}{4},\frac{1}{2}\right)$ for all $r \in \YY$. Then, the kernels $\mathfrak{f}$ defined in \eqref{equality_f^(u,l)_t,(m,s)} belong to $\mathcal{H}^2$.
\end{lemma}

\begin{proof}
    We must check the three conditions for a kernel to be admissible. The first condition, i.e., measurability, can be seen from the form of the kernels $\mathfrak{f}^{(r,s)}$. The second condition is verified since $\YY^2$ is a separable metric space in the Borel topology; see Remark \ref{rmk:meas}(1). The last condition follows from Lemma \ref{eq:finL2} below.
\end{proof}

\subsection{Rosenblatt distribution with sample paths in general $\mathbb{H} \otimes \Hi$}
We say that the $\Hi \otimes \Hi$-valued random variable $\mathfrak{Z}_U$ defined as
\begin{equation} \label{eq:Ztilde}
\mathfrak{Z}_U \doteq (U^* \otimes U^*) \mathcal{I}_{2}^{\sigma_U}(\mathfrak{f}),
\end{equation}
follows a Rosenblatt distribution when $\mathfrak{f}$ is as in \eqref{equality_f^(u,l)_t,(m,s)}, $U$ is a unitary operator as in the spectral decomposition \eqref{eq:UNU=D}, and $\sigma_U$ 
is defined by
\begin{equation*}
  \sigma_U(r,s) \doteq \EE\left( ((U \varepsilon)(r)) ( (U \varepsilon)(s) ) \right), \quad r, s \in \YY.
\end{equation*}
Note that in \eqref{eq:Ztilde}, we slightly abused notation, by viewing the sample paths of $\mathcal{I}_2^{\sigma_U}(\mathfrak{f})$ in the isomorphic space $L^2(\YY^2) \cong L^2(\YY) \otimes L^2(\YY)$.

\section{Proofs of Main Results}
\label{s:proofs}

In this section, we give the proofs of our main results. 
Theorems \ref{thm1:srd}, \ref{le:limit_theorem_L2}, and \ref{th:main_result} are respectively proven in 
Sections \ref{subs:lemma0}, \ref{subs:lemma}, and \ref{subs:proofmainresult}.

\subsection{Proof of Theorem \ref{thm1:srd}} \label{subs:lemma0}

The result can be inferred by following the proof of Theorem 5 in \cite{Mas}. We show that the arguments in \cite{Mas} remain true after replacing Assumption H.2 of \cite{Mas} by the weaker assumption \eqref{cond:srd-improved_1}. We note that \cite{Mas} considers two-sided (non-causal) series instead of one sided (causal) series, but this does not affect the analysis.

To complete the proof, we must show that Theorem 5 of \cite{Mas} holds under the weaker assumption \eqref{cond:srd-improved_1}. In turn, it suffices to prove that Lemma 8 in \cite{Mas} holds, i.e., we must check that
\begin{enumerate}[label=(\roman*)]
    \item $\sum\limits_{h = 0}^{\infty} \|\Gamma_{h+p -q} T \Gamma_h \|_{\operatorname{Tr}} < \infty$,
    \label{eq:estimate1}
    \item $\sum\limits_{h = 0}^{\infty} \|\Gamma_{h+q} T \Gamma_{h-p} \|_{\operatorname{Tr}} < \infty$, and
    \label{eq:estimate2}
    \item $\sum\limits_{h = 0}^{\infty} \| u_h (\Lambda - \Phi) u^*_{h+q} \|_{\operatorname{Tr}} < \infty$, \label{eq:estimate3}
\end{enumerate}
where $\Lambda, \Phi, \Gamma$ were defined in \eqref{def:def:Lambda-Phi}--\eqref{def:def:Aq-C}. We start with the third estimate \ref{eq:estimate3}. Note that
\begin{equation} \label{eq:999}
\begin{split}
   \sum_{h = 0}^{\infty} \| u_h (\Lambda - \Phi) u^{*}_{h+q} \|_{\operatorname{Tr}} 
   &\le 
   \left(\|\Lambda\|_{\operatorname{Tr}} + \|\Phi\|_{\operatorname{Tr}}\right) \sum_{h = 0}^{\infty} \|u_h\|_{\operatorname{op}} \|u_{h+q}\|_{\operatorname{op}} \\
   &\le 
   \left( \EE \|\varepsilon_0\|_{\mathbb{H}}^4 + 3\EE \|\varepsilon_0\|_{\mathbb{H}}^4 \right) \sum_{h = 0}^{\infty} \|u_h\|_{\operatorname{op}}^2 < \infty.
\end{split}
\end{equation}
For the second line in \eqref{eq:999} we used that, without loss of generality, we can select $\{u_j\}_{j \in \NN_0}$ in the representation \eqref{eq:linear_process} such that $\{\|u_j\|_{\operatorname{op}}\}_{j \in \NN_0}$ is decreasing. More precisely, we can define a renumeration $\{\tilde u_j\}_{j \in \NN_0}$ of $\{u_j\}_{j \in \NN_0}$ such that $\{ \tilde u_j \}_{j \in \NN_0}$ is decreasing in the operator norm. Then, take  $\tilde X_n \doteq  \sum_{j=0}^\infty \tilde u_j \varepsilon_{n-j}$. Note that, for all $n$, $X_n = \tilde X_n$ in law, since $\{\varepsilon_n\}_{n \in \ZZ}$ is i.i.d., and so we can work with $\tilde X_n$ instead of $X_n$ (and the same for their respective autocovariance operators).

We now turn to the first estimate \ref{eq:estimate1}, and \ref{eq:estimate2} follows by identical calculations. 
By denoting $M \doteq \sum_{j = 0}^{\infty} \|u_j\|_{\operatorname{op}}^{4/3} < \infty$ and since $\{ \|u_j \|_{\operatorname{op}}\}_{j \in \NN_0}$ is decreasing in $j$, we write
\begin{equation}
\begin{split}
    \sum_{h=0}^{\infty} \|\Gamma_{h+p -q} T \Gamma_h \|_{\operatorname{Tr}} 
    &\le \EE \|\varepsilon_0 \|_{\mathbb{H}}^4 \sum_{i = 0}^\infty \sum_{j = p}^\infty \sum_{h = 0}^\infty \|u_{h+i}\|_{\operatorname{op}}  \|u_{i}\|_{\operatorname{op}} \|u_{h+q+j}\|_{\operatorname{op}} \|u_{p+j}\|_{\operatorname{op}} \\
    &\le  \EE \|\varepsilon_0 \|_{\mathbb{H}}^4 \sum_{j} \sum_h  \|u_{h+q+j}\|_{\operatorname{op}} \|u_{p+j}\|_{\operatorname{op}} \|u_{h}\|_{\operatorname{op}}^{2/3} \sum_i \|u_{h+i}\|_{\operatorname{op}}^{1/3}  \|u_{i}\|_{\operatorname{op}} \\
    &\le M \EE \|\varepsilon_0 \|_{\mathbb{H}}^4 \sum_{h} \|u_{h}\|_{\operatorname{op}}^{2/3} \|u_{h+q+j}\|_{\operatorname{op}}^{2/3} \sum_j  \|u_{h+q+j}\|_{\operatorname{op}}^{1/3} \|u_{p+j}\|_{\operatorname{op}} \\
    &\le M^2 \EE \|\varepsilon_0 \|_{\mathbb{H}}^4 \sum_{h} \|u_{h}\|_{\operatorname{op}}^{2/3} \|u_{h+q+j}\|_{\operatorname{op}}^{2/3} \\
    &\le M^3 \EE \|\varepsilon_0 \|_{\mathbb{H}}^4.
\end{split}
\end{equation}
The estimate in the first line follows from calculations in pp. 127-128 of \cite{Mas}. Coupled with the strategy followed there, these estimates finish the proof.

\subsection{Proof of Theorem \ref{le:limit_theorem_L2}}
\label{subs:lemma}

In the proof of Theorem \ref{le:limit_theorem_L2}, we are concerned with the case $d(s) \in \left(\frac{1}{4},\frac{1}{2}\right), s \in \YY$. The proof is structured as follows. We first separate the sample autocovariances into their diagonal and off-diagonal terms. Lemma \ref{le:diagonal-1} below shows that only the off-diagonal terms contribute to the limit in the sought convergence of the sample autocovariance operators.
We then leverage Lemma \ref{Lemma_off_diagonals} below to identify the scaling limit (in $L^2(\mathbb{Y}^2 : \RR^{(H+1)})$) of the off-diagonal terms. Lemma \ref{Lemma_off_diagonals} is the key result in our analysis and its application gives the sought convergence result.
\par
The sample autocovariances can be separated into diagonal and off-diagonal parts as follows. For $h=0,\dots,H$, we have
\begin{equation} \label{eq:decomp}
\begin{aligned}
\widehat{\gamma}_{N,h}(r,s)-\gamma_{h}(r,s) 
&=	
\frac{1}{N} \sum_{n=1}^{N} 
X_{n+h}(r) X_{n}(s)-\sum_{j=0}^{\infty} (j+h+1)^{d(r)-1} (j+1)^{d(s)-1} \sigma(r,s)	
\\ & =
\frac{1}{N} \sum_{n=1}^{N} \sum_{j=0}^{\infty} u_{j+h}(r) u_{j}(s) (\varepsilon_{n-j}(r) \varepsilon_{n-j}(s)-\sigma(r,s))
\\ & \hspace{1cm} +
\frac{1}{N} \sum_{n=1}^{N} \sum_{j \neq i +h} u_{j}(r) u_{i}(s) \varepsilon_{n+h-j}(r) \varepsilon_{n-i}(s) 
\\ & = 
D_{N,h}(r,s)+O_{N,h}(r,s)
\end{aligned}
\end{equation}
with $u_{j}(r) = (j+1)^{d(r)-1}$ as in \eqref{eq:XninL2} and
\begin{equation} \label{eq:off_on_diag}
\begin{split}
    D_{N,h}(r,s) &\doteq \frac{1}{N} \sum_{n=1}^{N} \sum_{j=0}^{\infty} u_{j+h}(r) u_{j}(s) (\varepsilon_{n-j}(r) \varepsilon_{n-j}(s)-\sigma(r,s)), \\
    O_{N,h}(r,s) &\doteq \frac{1}{N} \sum_{n=1}^{N} \sum_{\substack{i,j=0,\dots, \infty \\ j \neq i + h}} u_{j}(r) u_{i}(s) \varepsilon_{n+h-j}(r) \varepsilon_{n-i}(s) .
\end{split}
\end{equation}
By Theorem 2.3 in \cite{Bosq2000:Linear}, it suffices to show
\begin{align} \label{eq:B-p1}
    \| \Xi_N (D_{N,h}, h = 0, \dots, H) \|_{L^2(\mathbb{Y}^2:\RR^{(H+1)}) } \overset{\P}{\to} 0,
\end{align}
\begin{align} \label{eq:B-p2}
    \Xi_N (O_{N,h}, h = 0, \dots, H)
\overset{\distr}{\to} (\mathfrak{R}, h = 0, \dots, H).
\end{align}
The convergences \eqref{eq:B-p1} and \eqref{eq:B-p2} follow by Lemmas \ref{le:diagonal-1} and \ref{le:diagonal-2} respectively.

\begin{lemma} \label{le:diagonal-1}
Let $\{X_{n}\}_{n \in \ZZ}$ be as in Theorem \ref{le:limit_theorem_L2}, with $d(s) \in \left(\frac{1}{4},\frac{1}{2}\right)$ for $\mu-$a.s. $s \in \mathbb{Y}$. Then, \eqref{eq:B-p1} holds.
\end{lemma}

\begin{proof}
By Markov's inequality, for all $\varepsilon > 0$, 
\begin{multline*}
    \P \left( \| \Xi_N  (D_{N,h}, h = 0, \dots, H) \|_{L^2(\mathbb{Y}^2 : \RR^{(H+1)})}  > \varepsilon \right) \leq \frac{1}{\varepsilon^2}
    \EE \| \Xi_N (D_{N,h}, h = 0, \dots, H) \|^2_{L^2(\mathbb{Y}^2 : \RR^{(H+1)}) } 
    \\= \frac{1}{\varepsilon^2}
    \int \int
    \EE \| N^{1-d(r,s)} (D_{N,h}(r,s), h = 0, \dots, H) \|^2_{\RR^{H+1}} \mu(dr)\mu(ds)
    \\= \frac{1}{\varepsilon^2}
    \sum_{h=0}^H 
    \int \int
    \EE | N^{1-d(r,s)} D_{N, h}(r, s) |^2 \mu(dr)\mu(ds).
\end{multline*}
It now suffices to show that the last term converges to $0$ as $N \to \infty$. Recalling \eqref{eq:off_on_diag}, note that we can rewrite the diagonal term as
\begin{equation}
D_{N, h} = \frac{1}{N} \sum_{n=1}^{N} Y_n \in L^2(\mathbb{Y}^2), \quad     Y_{n}(r,s) \doteq \sum_{j=0}^\infty \alpha^{(h)}_j [ \zeta_{j-n}] (r,s),
\end{equation}
where the sequences $\alpha_j^{(h)}, \zeta_j$ are respectively $L(L^2(\mathbb{Y}^2))$- and $L^2(\mathbb{Y}^2)$-valued elements given by 
\begin{equation}
    \alpha^{(h)}_j[f](r,s) \doteq u_{j+h}(r) u_{j}(s)f(r,s), \quad \zeta_{j}(r,s)\doteq \varepsilon_{j}(r) \varepsilon_{j}(s)-\sigma(r,s), \quad f \in L^2(\mathbb{Y}^2),\;  r , s \in \YY, \; j \ge 0.
\end{equation}
Moreover, note that $\{\zeta_j\}_{j \in \ZZ}$ are centered and i.i.d. with
\begin{equation} \label{eq:sigma-zeta}
    \EE \left( \zeta_{i}(r, s) \zeta_{j} (r, s)
    \right) =
    \begin{cases}
    \EE \varepsilon_0^2(r) \varepsilon_0^2(s) - (\sigma(r, s))^2 & i = j\\
    0 & i \neq j.
    \end{cases},
    \quad \sigma^2_\zeta(r,s) \doteq \EE \left( \zeta_{0}(r, s) \zeta_{0} (r, s)
    \right).
\end{equation}

We claim that there exists a constant $c>0$, such that
\begin{equation} \label{al:MI-all}
\begin{aligned}
    &\int \int
    \EE | N^{1-d(r,s)} D_{N, h}(r, s) |^2 \mu(dr,ds) = \int
     N^{2-2d(r,s)} \EE |D_{N, h}(r, s) |^2 \mu(dr,ds)
    \\
    &\leq 
    c \max\{ \log^2(N)/2, \log(N) \}
    \int \int N^{1-2d(r,s)} | \sigma_{\zeta}(r,s) |^2 \mu(dr) \mu(ds),
\end{aligned}
\end{equation}
where $\sigma_{\zeta}$ was defined in \eqref{eq:sigma-zeta}. Then, the quantity in the second line of \eqref{al:MI-all} converges to zero by the dominated convergence theorem and the observation that $\max\{ \log^2(N)/2, \log(N) \} N^{1-2d(r,s)} \to 0$ as $N \to \infty$ for a.s. $r,s \in \YY$, since $\log(N), \log^2(N)$ are slowly varying functions and $2d(r,s) < 1$ for $\mu$-a.s. $r,s \in \YY$. 

To show the inequality in \eqref{al:MI-all}, first note that
\begin{equation} \label{eq:456}
    \EE \left( D_{N, h}(r, s) \right)^2
    =
    \frac{1}{N} \EE \left( Y_{0}(r, s) \right)^2
    +
    \frac{2}{N^2} \sum_{n=1}^{N-1} \sum_{l=n+1}^N
    \EE \left( Y_{n}(r, s) Y_{l}(r, s) \right).
\end{equation}
We can see that
\begin{equation} \label{eq:0078}
\begin{split}
    \EE \left( Y_{0}(r, s) \right)^2 &= \sum_{i = 0}^\infty (i+1)^{d(r,r)-2}(i+h+1)^{d(s,s)-2} \sigma^2_{\zeta}(r,s) \\
    &\le \sigma^2_{\zeta}(r,s) \sum_{i = 0}^\infty (i+1)^{d(r,r)+d(s,s)-4} < \infty.
\end{split}
\end{equation}
Moreover, by using \eqref{eq:sigma-zeta}, it follows that
\begin{align}
    &\EE \left( Y_{0}(r, s) Y_{l}(r, s) \right) \nonumber
    \\&=
    \EE \left( 
    \sum_{i, j=0}^\infty \alpha^{(h)}_i[\zeta_{i}](r, s) \alpha^{(h)}_j[\zeta_{j-l}](r, s) 
    \right) \nonumber
    \\&= 
    \sum_{i, j=0}^\infty (i+1)^{d(r)-1} (i+h+1)^{d(s)-1} (j+1)^{d(r)-1} (j+h+1)^{d(s)-1} 
    \EE \left( \zeta_{i}(r, s) \zeta_{j-l} (r, s)
    \right) \nonumber
    \\&= 
    \sum_{j=0}^\infty (j+l+1)^{d(r)-1} (j+l+h+1)^{d(s)-1} (j+1)^{d(r)-1} (j+h+1)^{d(s)-1} 
    \sigma^2_{\zeta}(r, s) \nonumber
    \\&\leq 
    \sum_{j=0}^\infty (j+l+1)^{d(r, s)-2} (j+1)^{d(r, s)-2}
    | \sigma_{\zeta}(r, s) |^2 
    \leq | \sigma_{\zeta}(r, s) |^2
    \sum_{j=0}^\infty (j+l+1)^{-1} (j+1)^{-1}
    . \label{eq:987}
\end{align}
Moreover, by the proof of Proposition 1 in \cite{Chara_Rack_Central},
\begin{equation*}
    \sum_{j=0}^\infty (j+l+1)^{-1} (j+1)^{-1} 
    \leq
    (l+1)^{-1} + l^{-1} \left(
    \log \left( \frac{l+1}{2} \right) + \int_1^{\infty} ( x(x+1))^{-1} dx
    \right) \\
    \leq c l^{-1} \log(l),
\end{equation*}
so that, there is a constant, $c>0$, with
\begin{equation} \label{eq:988}
    \sum_{j=0}^\infty (j+l+1)^{-1} (j+1)^{-1} 
    \leq c l^{-1} \log(l).
\end{equation}
Then the last term in \eqref{eq:456} can be estimated upon noticing that
\begin{align}
    \frac{1}{N^2} \sum_{n=1}^{N-1} \sum_{l=n+1}^N
    \EE \left( Y_{n}(r, s) Y_{l}(r, s) \right)
    &=
    \frac{1}{N^2} \sum_{l=1}^{N-1} (N-l) 
    \EE \left( Y_{0}(r, s) Y_{l}(r, s) \right) \nonumber
    \\&\leq c
    \frac{1}{N^2} \sum_{l=1}^{N-1} (N-l)
    l^{-1} \log(l)
    | \sigma_{\zeta}(r, s) |^2
        \label{eq:989}
    \\&\leq
    c
    \frac{1}{N} \sum_{l=1}^{N-1} (1-\frac{l}{N})
    l^{-1} \log(l)
    | \sigma_{\zeta}(r, s) |^2
    \nonumber \\
    &\leq c | \sigma_{\zeta}(r, s) |^2 \frac{\log N}{N} \sum_{l=1}^{N-1} l^{-1} \nonumber
    \\&\leq c
    \frac{1}{N} \max\{ \log^2(N)/2, \log(N) \} | \sigma_{\zeta}(r, s) |^2,
    \label{eq:990}
\end{align}
where \eqref{eq:989} follows from \eqref{eq:987}--\eqref{eq:988}. Then, \eqref{al:MI-all} follows from \eqref{eq:456}, \eqref{eq:0078}, and \eqref{eq:990}.
\end{proof}

\begin{remark}
In the case $\esssup_{s \in \YY} d(s) < \frac{1}{2}$, $\{\alpha_j^{(h)}\}_{j \in \NN}$ in the proof of Lemma \ref{le:diagonal-1} are absolutely summable in the operator norm. Then, Lemma \ref{le:diagonal-1}, as well as the asymptotic normality for $D_{N,h}$, follows from Theorem 2 in \citet{Merlevede1997}. This argument is no longer available when $\esssup_{s \in \YY} d(s) = \frac{1}{2}$, thus justifying the proof strategy pursued here.
\end{remark}

\begin{lemma} \label{le:diagonal-2}
Let $\{X_{n}\}_{n \in \ZZ}$ be as in Theorem \ref{le:limit_theorem_L2}, with $d(s) \in \left(\frac{1}{4},\frac{1}{2}\right)$ for $\mu-$a.s. $s \in \mathbb{Y}$. Then, \eqref{eq:B-p2} holds.
\end{lemma}

\begin{proof}
To investigate the asymptotic behavior of the off-diagonal terms in \eqref{eq:B-p2}, we write
\begin{align*}
N^{1-d(r,s)} O_{N,h}(r,s)
& =
N^{-d(r,s)} \sum_{n=1}^{N} \sum_{\substack{i,j = 0 \\ j \neq i+h}}^\infty u_{j}(r) u_{i}(s) \varepsilon_{n-(j-h)}(r) \varepsilon_{n-i}(s)\\
& =
N^{-d(r,s)} \sum_{n=1}^{N} \sum_{\substack{j_1 = -h, j_2 = 0 \\ j_{1} \neq j_{2}}}^\infty
u_{j_{1}+h}(r) u_{j_{2}}(s) \varepsilon_{n-j_{1}}(r) \varepsilon_{n-j_{2}}(s)	\\
& =
N^{-d(r,s)} \sum_{n=1}^{N} \sum_{\substack{j_1 = -\infty, \dots, n+h  \\ j_2 = -\infty, \dots, n \\ j_{1} \neq j_{2}}}  u_{n+h-j_{1}}(r) u_{n-j_{2}}(s) \varepsilon_{j_{1}}(r) \varepsilon_{j_{2}}(s)	\\
& = \sum_{\substack{j_1,j_2 = -\infty \\ j_{1} \neq j_{2}}}^\infty N^{-d(r,s)} \sum_{n=1}^{N} u_{n+h-j_{1}}(r) u_{n-j_{2}}(s) \varepsilon_{j_{1}}(r) \varepsilon_{j_{2}}(s)  \mathds{1}_{\{n+h-j_{1} \geq 0 \}} \mathds{1}_{\{n-j_{2} \geq 0 \}}	\\
& =
\sum_{\substack{j_{1} \neq j_{2}}} C_{N,h}(j_{1},j_{2},r,s)
\varepsilon_{j_{1}}(r) \varepsilon_{j_{2}}(s),
\end{align*}
where we have used the changes of variables $j_1 = j - h$ ($j_2= i$), and $j_k = n - j_k$, $k=1,2$, for the second and third lines respectively, and where
\begin{equation}	\label{eq:Cn}
C_{N,h}(j_{1},j_{2},r,s) 
\doteq
N^{-d(r,s)} \sum_{n=1}^{N}  u_{n+h-j_{1}}(r) u_{n-j_{2}}(s) \mathds{1}_{\{n+h-j_{1} \geq 0 \}} \mathds{1}_{\{n-j_{2} \geq 0 \}}.
\end{equation}
Then,
\begin{equation}
    \Xi_N (O_{N,h}, h=0,\dots, H)(r,s) = 
    \sum_{\substack{j_{1} \neq j_{2}}} C_{N}(j_{1},j_{2},r,s) \varepsilon_{j_{1}}(r) \varepsilon_{j_{2}}(s), \; r,s \in \YY,
\end{equation}
with $C_{N}(j_{1},j_{2},\cdot,\cdot) \in L^2(\mathbb{Y}^2:\RR^{(H+1)})$ defined by
\begin{equation} \label{eq:CN-vector}
    C_{N}(j_{1},j_{2},r,s) = ( C_{N,0}(j_{1},j_{2},r,s) , \dots, C_{N,H}(j_{1},j_{2},r,s) )' .
\end{equation}
Moreover, define for $r,s \in \mathbb{Y}$, $h=0,\dots,H$, and $x_1,x_2 \in \RR$,
\begin{equation} \label{eq:def-cntilde}
\widetilde{C}_{N}(x_{1},x_{2},\cdot,\cdot) \doteq N C_{N}(\lceil x_{1}N \rceil,\lceil x_{2}N \rceil ,\cdot,\cdot), \quad \widetilde{C}_{N}^\sigma(x_1,x_2,r,s) \doteq \widetilde{C}_{N}(x_1,x_2,r,s) \sigma(r) \sigma(s).
\end{equation}
Then, from \eqref{eq:Cn},
\begin{align}
&
\widetilde C_{N,h}(x_1,x_2,r,s) \nonumber
\\& =
N^{1-d(r,s)} \sum_{n=1}^{N}  u_{n+h-\lceil x_{1}N \rceil}(r)  u_{n -\lceil x_{2}N\rceil }(s) 
\mathds{1}_{\{n+h-\lceil x_{1}N \rceil \geq 0 \}} \mathds{1}_{\{n-\lceil x_{2}N \rceil \geq 0 \}} \nonumber \\
& =
N^{2-d(r,s)} 
\int_{0}^{1} u_{ \lceil vN \rceil + h -\lceil x_{1}N \rceil}(r) u_{\lceil vN \rceil-\lceil x_{2}N \rceil}(s) 
\mathds{1}_{\{\lceil vN \rceil+h-\lceil x_{1}N \rceil \geq 0 \}} \mathds{1}_{\{\lceil vN \rceil-\lceil x_{2}N \rceil \geq 0 \}}
dv. \label{eq:659}
\end{align}
We seek to apply Lemma \ref{Lemma_off_diagonals} to identify the weak limits of $\Xi_N (O_{N,h}, h=0,\dots, H)$. For this reason, we set 
\begin{equation}
    \widebar{\mathfrak{f}}^\sigma \doteq ( \mathfrak{f}^\sigma, h=0,\dots, H)', \quad \mathfrak{f}^\sigma(x,y,r,s) \doteq \mathfrak{f}^{(r,s)}(x,y) \sigma(r) \sigma(s)
\end{equation}
as in \eqref{eq:defn-f-hat}, where $\mathfrak{f}$ was defined in \eqref{equality_f^(u,l)_t,(m,s)}. Here we maintain the notation $\sigma$ to remember the dependence on the correlation structure of the underlying noise. Then, the previous relation implies that,
\begin{equation} \label{eq:cond-for-lemma}
\begin{split}
&\left\| \widetilde{C}^\sigma_{N} - \widebar{\mathfrak{f}}^\sigma \right\|^2_{L^{2}(\RR^2: L^2(\mathbb{Y}^2:\RR^{(H+1)}))}	 
\\&=
\int_{\RR^2} \int_{\mathbb{Y}^2} 	
\Bigg\| \Bigg( \int_{0}^{1} N^{2-d(r,s)} 
u_{\lceil vN \rceil+h- \lceil x_{1}N]}(r) u_{\lceil vN \rceil-\lceil x_{2}N \rceil}(s) 
\mathds{1}_{\{\lceil vN \rceil +h-\lceil x_{1}N \rceil \geq 0 \}} \mathds{1}_{\{\lceil vN \rceil-\lceil x_{2}N \rceil \geq 0 \}}
\\
&\hspace{2cm} - 
(v-x_1)_+^{d(r)-1} (v-x_2)_+^{d(s)-1}dv, h =0,\dots,H \Bigg) \Bigg\|_{\RR^{H+1}}^2 \sigma^2(r) \sigma^2(s) dr ds dx_1 dx_2.
\end{split}
\end{equation}
Note that to show $\left\| \widetilde{C}^\sigma_{N} - \widebar{\mathfrak{f}}^\sigma \right\|^2 \to 0$, it remains to show that the integrand is (i) integrable, and (ii) converges pointwise (in $r,s,x_1,x_2$, as $N \to \infty$) to 0. 
Upon establishing these two claims, the DCT will ensure that the integral tends to zero. 

For (i), recall that from Lemma \ref{eq:finL2} we have that $\mathfrak{f}^\sigma \in {L^{2}(\RR^2: {L^2(\mathbb{Y}^2)})}$. For the sake of notational simplicity, we set $H = 0$ (and so $h = 0$), and we show that $\widetilde C_N^\sigma \in {L^{2}(\RR^2: {L^2(\mathbb{Y}^2)})}$. The case $h \in \NN$ follows from similar calculations. First, note that 
\begin{equation} \label{eq:677}
    \begin{aligned}
        &N^{1-d(r)} u_{\lceil vN\rceil -\lceil x_{1}N\rceil}(r) \mathds{1}_{\{\lceil vN \rceil-\lceil x_{1}N \rceil \geq 0 \}} \\
        &\hspace{1cm}= \left(\frac{1}{N}\right)^{d(r)-1} \mathds{1}_{\{\lceil vN \rceil = \lceil x_{1}N \rceil  \}} + 
        \left( \frac{\lceil vN \rceil}{N} - \frac{\lceil x_1 N \rceil}{N}  + \frac{1}{N}    \right)^{d(r)-1} \mathds{1}_{\{\lceil vN \rceil -\lceil x_{1}N \rceil \geq 1 \}}, 
    \end{aligned}
\end{equation}
where it follows that, for all $x_1 \in \RR, r \in \YY, N \in \NN$,
\begin{equation} \label{eq:678}
        \left( \frac{\lceil vN \rceil}{N} - \frac{\lceil x_1 N \rceil}{N}  + \frac{1}{N}    \right)^{d(r)-1} \mathds{1}_{\{\lceil vN \rceil-\lceil x_{1}N \rceil \geq 1 \}}
        \le \left( v  - x_1   \right)_+^{d(r)-1},
\end{equation}
by recalling that $v \le \frac{\lceil vN \rceil}{N}$, $\frac{\lceil x_1 N\rceil}{N} \le x_1 + \frac{1}{N}$, and that $d(r) - 1 < 0$. Analogous calculations show that
\begin{equation} \label{eq:679}
    \begin{aligned}
    &
    N^{1-d(s)} u_{\lceil vN \rceil -\lceil x_{2}N \rceil}(s) \mathds{1}_{\{\lceil vN \rceil-\lceil x_{2}N \rceil \geq 0 \}}  \\
    &\hspace{1cm}=
    \left(\frac{1}{N}\right)^{d(s)-1} \mathds{1}_{\{\lceil vN \rceil = \lceil x_{2}N \rceil  \}} + 
        \left( \frac{\lceil vN \rceil}{N} - \frac{\lceil x_2 N \rceil}{N}  + \frac{1}{N}    \right)^{d(s)-1} \mathds{1}_{\{\lceil vN \rceil -\lceil x_{1}N \rceil \geq 1 \}} \\
    &\hspace{1cm}\le \left(\frac{1}{N}\right)^{d(s)-1} \mathds{1}_{\{\lceil vN \rceil = \lceil x_{2}N \rceil  \}} +  \left( v  - x_2   \right)_+^{d(s)-1}.
    \end{aligned}
\end{equation}
Combining \eqref{eq:659}, \eqref{eq:677}, \eqref{eq:678}, and \eqref{eq:679}, we see that 
\begin{equation*} 
    \left\| \widetilde C^\sigma_{N} \right\|^2_{L^{2}(\RR^2: L^2(\mathbb{Y}^2))} 
    \le 4 \left(  \mathcal{R}_{1,N} + \mathcal{R}_{2,N} + \mathcal{R}_{3,N} +
    \left\| \mathfrak{f}^\sigma \right\|^2_{L^{2}(\RR^2: L^2(\mathbb{Y}^2))} \right),
\end{equation*}
where
\begin{equation}
\begin{split}
    \mathcal{R}_{1,N} &\doteq \left\|\sigma(r) \sigma(s)  N^{2-d(r,s)} \int_0^1 \mathds{1}_{\{\lceil vN \rceil = \lceil x_{1}N \rceil  \}} \mathds{1}_{\{\lceil vN \rceil = \lceil x_{2}N \rceil   \}} dv \right\|^2_{L^{2}(\RR^2: L^2(\mathbb{Y}^2))}, \\
    \mathcal{R}_{2,N} &\doteq \left\|\sigma(r) \sigma(s) N^{1-d(r)} \int_0^1 (v-x_2)_+^{d(s)-1} \mathds{1}_{\{\lceil vN \rceil = \lceil x_{1}N \rceil  \}}  dv \right\|^2_{L^{2}(\RR^2: L^2(\mathbb{Y}^2))}, \\
    \mathcal{R}_{3,N} &\doteq \left\| \sigma(r) \sigma(s) N^{1-d(s)} \int_0^1 (v-x_1)_+^{d(r)-1} \mathds{1}_{\{\lceil vN \rceil = \lceil x_{2}N \rceil  \}}  dv 
 \right\|^2_{L^{2}(\RR^2: L^2(\mathbb{Y}^2))}.
\end{split}
\end{equation}
We now show that $\mathcal{R}_{i,N} \to 0$ as $N \to \infty$ for $i=1,2,3$. We first investigate $\mathcal{R}_{1,N}$. Note that
\begin{equation} \label{eq:790}
\begin{split}
    \mathcal{R}_{1,N} &= \int_{\YY^2} \sigma^2(r) \sigma^2(s) N^{4-2d(r,s)} \int_{\RR^2} \int_0^1 \int_0^1 \mathds{1}_{\{\lceil v_1 N \rceil = \lceil x_{1}N \rceil, \lceil v_1N \rceil = \lceil x_{2}N \rceil , \lceil v_2N \rceil = \lceil x_{1}N \rceil , \lceil v_2 N \rceil = \lceil x_{2}N \rceil  \}} \\
    &\hspace{7cm} \times dv_1 dv_2 dx_1 dx_2 \mu(dr) \mu(ds) \\
    &= \int_{\YY^2} \sigma^2(r) \sigma^2(s) N^{1-2d(r,s)} \mu(dr) \mu(ds) \\
    &= \left( \int_{\YY} \sigma^2(r) N^{1/2-2d(r)} \mu(dr) \right)^2,
\end{split}
\end{equation}
where the second equality follows upon noticing that
\begin{multline*}
     \int_{\RR^2} \int_0^1 \int_0^1 \mathds{1}_{\{\lceil v_1 N \rceil = \lceil x_{1}N \rceil, \lceil v_1N \rceil = \lceil x_{2}N \rceil , \lceil v_2N \rceil = \lceil x_{1}N \rceil , \lceil v_2 N \rceil = \lceil x_{2}N \rceil  \}} \\
     = \sum_{i=0}^{N-1} \int_{i/N}^{(i+1)/N} \int_{i/N}^{(i+1)/N} \int_{i/N}^{(i+1)/N} \int_{i/N}^{(i+1)/N} dv_1 dv_2 dx_1 dx_2 =
     N^{-3}.
\end{multline*}
Since $d(r) > \frac{1}{4}$, it follows from \eqref{eq:790} and DCT that $ \mathcal{R}_{1,N} \to 0$. We now turn to proving $\mathcal{R}_{2,N} \to 0$ and the calculations are analogous (but symmetrical in $r$ and $s$) for $\mathcal{R}_{3,N} \to 0$. First, observe that
\begin{multline} \label{eq:455}
     \left( N^{1-d(r)} \int_0^1 (v-x_2)_+^{d(s)-1} \mathds{1}_{\{\lceil vN \rceil = \lceil x_{1}N \rceil  \}}  \right)^2 \\
     = N^{2-2d(r)} \int_0^1 \int_0^1 (v_1-x_2)_+^{d(s)-1} (v_2-x_2)_+^{d(s)-1}  \mathds{1}_{\{\lceil v_1N \rceil = \lceil x_{1}N \rceil  \}}  \mathds{1}_{\{\lceil v_2N \rceil = \lceil x_{1}N \rceil  \}} dv_1 dv_2,
\end{multline}
and moreover
\begin{equation} \label{eq:4565}
    \begin{aligned}
    &
    \int_{\RR}  \mathds{1}_{\{\lceil v_1N \rceil = \lceil x_{1}N \rceil  \}}  \mathds{1}_{\{\lceil v_2N \rceil = \lceil x_{1}N \rceil  \}} dx_1 
    \\&\le \mathds{1}_{\{\lceil v_1N \rceil = \lceil v_{2} N \rceil  \}} \sum_{i=0}^{N-1} \int_{i/N}^{(i+1)/N}  \mathds{1}_{\{v_1 \in \left(\frac{i}{N}, \frac{i+1}{N}\right) \}} dx_1
    = \mathds{1}_{\{\lceil v_1N \rceil = \lceil v_{2} N \rceil  \}} \frac{1}{N}.
    \end{aligned}
\end{equation} 
In addition, we see that
\begin{align} 
    &
    \int_0^1 \int_0^1 \mathds{1}_{\{\lceil v_1N \rceil = \lceil v_{2} N \rceil  \}} \int_\RR (v_1-x_2)_+^{d(s)-1} (v_2-x_2)_+^{d(s)-1} dx_2 dv_2 dv_1 
    \nonumber
    \\&
    = \int_0^1 \int_0^1 \mathds{1}_{\{\lceil v_1N \rceil = \lceil v_{2} N \rceil  \}} \int_{- \min\{v_1,v_2\}}^\infty (v_1+x_2)_+^{d(s)-1} (v_2+x_2)_+^{d(s)-1} dx_2 dv_2 dv_1 
    \nonumber
    \\&
    = \int_0^1 \Bigg[ \int_0^{v_1} \mathds{1}_{\{\lceil v_1N \rceil = \lceil v_{2} N \rceil  \}} \int_{- v_2}^\infty  (v_1+x_2)^{d(s)-1} (v_2+x_2)^{d(s)-1} dx_2 dv_2 
    \nonumber
    \\&\hspace{1cm}+  \int_{v_1}^1 \mathds{1}_{\{\lceil v_1 N \rceil = \lceil v_{2} N \rceil  \}} \int_{- v_1}^\infty  (v_1+x_2)^{d(s)-1} (v_2+x_2)^{d(s)-1} dx_2 dv_2  \Bigg] dv_1
    \nonumber\\&
    = \int_0^1 \Bigg[ \int_0^{v_1} \mathds{1}_{\{\lceil v_1 N \rceil = \lceil v_{2} N \rceil  \}} \int_{0}^\infty  z^{d(s)-1} (z+ (v_1-v_2))^{d(s)-1} dz dv_2 
    \nonumber
    \\&\hspace{1cm}+  \int_{v_1}^1 \mathds{1}_{\{\lceil v_1N \rceil = \lceil v_{2} N \rceil  \}} \int_0^\infty  z^{d(s)-1} (z+(v_2-v_1))^{d(s)-1} dx_2 dv_2  \Bigg] dv_1 
    \nonumber\\
     &=\int_0^1  \int_0^1 \mathds{1}_{\{\lceil v_1N \rceil = \lceil v_{2} N \rceil  \}} \int_{0}^\infty  z^{d(s)-1} (z+ |v_1-v_2|)^{d(s)-1} dz dv_2 dv_1   \nonumber\\
    &=  \int_0^1 \int_0^1 \mathds{1}_{\{\lceil v_1N \rceil = \lceil v_{2} N \rceil  \}} |v_2 - v_1|^{2 d(s)-1} B(d(s),1-2d(s)) dv_1 dv_2
    \nonumber
    \\&\le B(d(s),1-2d(s)) N^{-2d(s)}.
    \label{eq:457}
\end{align}
Here, the first equality followed from the change of variables $-x_2 \mapsto x_2$; the fifth and sixth lines follow by the change of variables $x_2 + v_2 = z$ and $x_2 + v_1 =z$ respectively; finally the last line follows upon noticing that, on the event $\mathds{1}_{\{\lceil v_1N \rceil = \lceil v_{2} N \rceil  \}}$, we have $|v_2 - v_1| \le N^{-1}$.

By combining \eqref{eq:455}, \eqref{eq:4565}, and \eqref{eq:457} we have that
\begin{equation}
\begin{split}
    \mathcal{R}_{2,N} &\le \int_{\YY^2} \sigma^2(r) \sigma^2(s) N^{1-2d(r)} \int_0^1 \int_0^1 \mathds{1}_{\{\lceil v_1N \rceil = \lceil v_{2} N \rceil  \}} \\
    &\hspace{4cm} \times \int_\RR (v_1-x_2)_+^{d(s)-1} (v_2-x_2)_+^{d(s)-1} dx_2 dv_1 dv_2 \mu(dr) \mu(ds) \\
    &\le \int_{\YY^2} \sigma^2(r) \sigma^2(s) N^{1-2d(r,s)} B(d(s), 1 - 2d(s)) \mu(dr) \mu(ds) \\
    &\le \int_{\YY} \sigma^2(r) N^{1/2 - 2d(r)} \mu(dr) \int_{\YY} \sigma^2(s)   N^{1/2 - 2d(s)} \frac{1}{d(s) \left(1-2d(s)\right)} \mu(ds) \\
    &\le 4\int_{\YY} \sigma^2(r) N^{1/2 - 2d(r)} \mu(dr) \int_{\YY}    N^{1/2 - 2d(s)} \frac{\sigma^2(s)}{1-2d(s)} \mu(ds),
\end{split}
\end{equation}
where the third inequality follows from the bound $0 \le B(x,y) \le \frac{1}{xy}$; see Theorem 1.1 in \cite{FroRat22}. Then, $\mathcal{R}_{2,N} \to 0$ as $N \to \infty$ from Assumption \eqref{eq:lemmaint2}.

The calculations above imply that
\begin{align}
\left\| \widetilde{C}^\sigma_{N} -
    \widebar{\mathfrak{f}}^\sigma \right\|^2_{L^{2}(\RR^2: L^2(\mathbb{Y}^2))} &\le 
    2\left\| \widetilde{C}^\sigma_{N}\right\|^2_{L^{2}(\RR^2: L^2(\mathbb{Y}^2))} 
    + 
    2\left\| \widebar{\mathfrak{f}}^\sigma \right\|^2_{L^{2}(\RR^2: L^2(\mathbb{Y}^2))}  \nonumber
    \\&\le 
    10\left\| \widebar{\mathfrak{f}}^\sigma \right\|^2_{L^{2}(\RR^2: L^2(\mathbb{Y}^2))} + 8\sum_{i=1}^3 \mathcal{R}_{i,N} < \infty. \label{eq:966}
\end{align}
We turn to point (ii). Note that from the definition of $\{u_j\}_{j \in \NN}$, for all $v,x_1,x_2 \in \RR$ and $r,s \in \YY$, as $N \to \infty$
\begin{equation}
\begin{split}
    N^{1-d(r)} u_{\lceil vN\rceil + h -\lceil x_{1}N\rceil}(r)  
    \mathds{1}_{\{n+h-j_{1} \geq 0 \}} \sigma(r) \to (v-x_1)_+^{d(r)-1} \sigma(r),
    \\
    N^{1-d(s)} u_{\lceil vN \rceil-\lceil x_{2}N \rceil}(s)  \mathds{1}_{\{n-j_{2} \geq 0 \}} \sigma(s) \to (v-x_2)_+^{d(s)-1} \sigma(s).
\end{split}
\end{equation}
The a.s. pointwise convergence (in $(r,s,x_1,x_2)$ of $\widetilde{C}_{N,h}(x_1,x_2,r,s)$ to $\mathfrak{f}^\sigma$) follows by DCT by using similar bounds as in the proof of (i), and noticing that these bounds are finite for $\lambda^{\otimes 2} \otimes \mu^{\otimes 2}$ a.s. $(x_1,x_2,r,s)$.

Combining items (i), (ii), and DCT, we can infer that 
\begin{equation}
\| \widetilde{C}_{N}^\sigma - \bar{\mathfrak{f}}^\sigma \|^2_{L^{2}(\RR^2: {L^2(\mathbb{Y}^2:\RR^{H+1})})}	 \to 0
\end{equation}
This shows that the conditions in Lemma \ref{Lemma_off_diagonals} are satisfied, finishing the proof of \eqref{eq:B-p2}.
\end{proof}

\subsection{Proof of Theorem \ref{th:main_result}}
\label{subs:proofmainresult}

The proof of Theorem \ref{th:main_result} is based on Theorem \ref{le:limit_theorem_L2}. We show that the $\Hi$-valued linear process \eqref{eq:linear_process} with \eqref{eq:linear_process2}
can be written as a continuous operator applied to an $L^2(\mathbb{Y})$-valued linear process \eqref{eq:linear_process}. Analogously, the sample autocovariance operator can be written as a continuous operator applied to the sample autocovariances of an $L^2(\mathbb{Y})$-valued linear process. Then, the continuous mapping theorem and Theorem \ref{le:limit_theorem_L2} give the desired result. 
\par
Recall the spectral theorem from \eqref{eq:UNU=D}. The self-adjoint operator $T$ can be decomposed into a multiplication operator $D_d$ and a unitary operator $U$. Then, 
\begin{equation} \label{eq:process_normal_perator_after_interchange}
\begin{split}
X_{n} =
\sum_{j=0}^{\infty} (j+1)^{T-I} \varepsilon_{n-j} &= \sum_{j=0}^\infty \exp\left( U^* (D_d - I) U \log(j+1)   \right) \varepsilon_{n-j} \\
&= \sum_{j=0}^\infty \sum_{\kappa = 0}^\infty U^* \frac{(D_d-I)^\kappa \log(j+1)^\kappa}{\kappa!} U \varepsilon_{n-j} \\
&= U^{*}\left( \sum_{j=0}^{\infty} (j+1)^{D_d-I}(U\varepsilon_{n-j})\right),
\end{split}
\end{equation}
where we used, for $T \in L(\Hi)$ and $\lambda>0$,   $\lambda^{T}=e^{T\log(\lambda)}$ and $e^{T}=\sum_{j=0}^{\infty} T^{j}/j!$ in the first and second line, respectively. Therefore,  
\begin{equation} \label{equality_process_Z}
X_n = U^* Z_n
\hspace{0.2cm}
\text{ with }
\hspace{0.2cm}
Z_{n} \doteq \sum_{j=0}^{\infty} (j+1)^{D_d-I}(U\varepsilon_{n-j}) \in L^2(\mathbb{Y}).
\end{equation}

Recall from Section \ref{subsec:prelim-linear} that the decomposition \eqref{equality_process_Z}, allows us to infer some properties of $\{X_{n}\}_{n \in \ZZ}$ based on $Z_n$.
In particular, if $Z_n$ satisfies \eqref{eq:integ-cond-well-def} with $\sigma^2(r) = \EE ( (U \varepsilon(r)) (U \varepsilon(r)))$, the unitarity of $U$ implies that $\{X_{n}\}_{n \in \ZZ}$ converges $\P$-almost surely.
\par
The interchangeability of the series and $U^*$ in \eqref{eq:process_normal_perator_after_interchange} is a consequence of the almost sure convergence of $Z_{n}$ in \eqref{equality_process_Z} and the boundedness of the unitary operator $U$. 
The process $\{Z_{n}\}_{ n \in \ZZ }$ satisfies the assumptions of Theorem \ref{le:limit_theorem_L2}: 
The sequence $\{ U\varepsilon_{j} \}_{j \in \ZZ}$ is an i.i.d.\ sequence with finite second and fourth moments
since $U: \Hi \to L^2(\mathbb{Y})$ is a unitary operator and $\{ \varepsilon_{j} \}_{j \in \ZZ}$ is assumed to be an i.i.d.\ sequence with finite second and fourth moments.
\par 
In view of \eqref{equality_process_Z}, we can rewrite the sample autocovariance operators of lag $h$ as
\begin{equation*}
\begin{aligned}
\widehat{\Gamma}_{N,h}-\Gamma_{h}
& =		
\frac{1}{N} \sum_{n=1}^{N} X_{n+h} \otimes  X_{n}-\EE(X_{h} \otimes X_{0}) 
\\ &=
\frac{1}{N} \sum_{n=1}^{N} ( (U^{*}Z_{n+h}) \otimes (U^{*}Z_{n})-E( (U^{*}Z_{h}) \otimes (U^{*}Z_{0}))) 
\\ &=
\frac{1}{N} \sum_{n=1}^{N}  \left( (U^* \otimes U^*) (Z_{n+h} \otimes Z_n) -(U^* \otimes U^*) \EE(Z_{h} \otimes Z_{0})\right) 
\\ &=
(U^* \otimes U^*) \left[ \frac{1}{N} \sum_{n=1}^{N}  \left(  Z_{n+h} \otimes Z_n - \EE(Z_{h} \otimes Z_{0})\right)  \right].
\end{aligned}
\end{equation*}
To interchange the summation and integration operations with $U$ in the calculations above, we used the fact that $U$ is a unitary, bounded, and linear operator. 
\par
Recall the scaling operators $\Xi_N, \Delta^U_N$ defined in \eqref{eq:Xi-def} and \eqref{eq:scaling:Delta}. It follows that the normalized sample autocovariance operator \eqref{eq:SampleAutoShort} of $X_{n}$ can be written as
\begin{equation} \label{eq:bla}
\begin{aligned}
\Delta^U_{N} (\widehat{\Gamma}_{N,h}-\Gamma_{h}) & =	(U^* \otimes U^*) 
 \left[ \Xi_N \left( \frac{1}{N} \sum_{n=1}^{N}  \left(  Z_{n+h} \otimes Z_n - \EE(Z_{h} \otimes Z_{0})\right)  \right)   \right]
\\ & \xrightarrow{d}	
(U^* \otimes U^*) \mathcal{I}^{\sigma_U}_2(f)
\end{aligned}
\end{equation}
with 
\begin{equation*}
  \sigma_U(r,s) \doteq \EE\left( ((U \varepsilon)(r)) ( (U \varepsilon)(s) ) \right), \quad r, s \in \YY.
\end{equation*}
The weak convergence
\begin{equation*}
    \Xi_N \left( \frac{1}{N} \sum_{n=1}^{N}  \left(  Z_{n+h} \otimes Z_n - \EE(Z_{h} \otimes Z_{0})\right)  \right)  \xrightarrow{d} \mathcal{I}^{\sigma_U}_2(f)
\end{equation*}
is in $L^2(\YY^2)$ and follows from an application of Theorem \ref{le:limit_theorem_L2}.
The weak convergence in \eqref{eq:bla} is in $\Hi \otimes \Hi$ and holds from the continuous mapping theorem, since $U^* \otimes U^*$ is bounded and linear. Note that the second line in \eqref{eq:bla} includes a slight abuse in notation, by identifying $L^2(\YY^2) \cong L^2(\YY) \otimes L^2(\YY)$.
\par
We can conclude that the continuous mapping theorem, \eqref{eq:bla} and Lemma \ref{le:limit_theorem_L2} give the desired convergence result.

\section{Some Open Questions}

The present work identifies the scaling limit of the autocovariance operator for a linear, discrete-time stochastic process taking values in a separable Hilbert space. Our focus has been on processes exhibiting long-range dependence. The tools developed herein exhibit considerable flexibility, and we conjecture that they could be adapted to identify scaling limits in a variety of other scenarios. For example:

\begin{openquest}
 In the finite-dimensional setting, it is possible to establish the convergence of higher-order statistics to Hermite processes of corresponding order; see Section 5.6 in \cite{PipirasTaqqu}. These processes can be represented as multiple Wiener-Itô stochastic integrals. We believe that the construction of the integrals in Section \ref{se:biW} can be extended to higher-order multiple Wiener-Itô integrals. This would pave the way for a full generalization of Theorem 5.6.3 in \cite{PipirasTaqqu} to the setting of linear processes taking values in $L^2$.
\end{openquest}

\begin{openquest}
In the second regime, we imposed a specific representation for the coefficients $u_j = (j + 1)^{D_d - I}$ in the process. This can easily be generalized to the case $u_j = \ell(j) (j + 1)^{D_d - I}$, where $\ell(j)$ is a slowly varying function. Is it possible to relax this assumption further?
\end{openquest}

On the other hand, our proof techniques have certain limitations. We highlight two such limitations below:

\begin{openquest}
The ``boundary'' case where $d(r) = \frac{1}{2}$ for $r \in A$ with $\mu(A) > 0$ has not been addressed. We conjecture that when $d(r) \in (0, 1/4]$, the limiting random variable remains Gaussian. However, by analogy with the finite-dimensional case, we expect that the appropriate scaling operator becomes logarithmic on the set $A$ where $d(r) = \frac{1}{4}$ for $r \in A$.
\end{openquest}

\begin{openquest}
In view of Corollary \ref{cor:srd-L^2} and Theorem \ref{le:limit_theorem_L2}, A natural question concerns the ``mixed'' case where $d(s) \in \left(0, \frac{1}{4}\right)$ for $s \in A$, $d(s) \in \left(\frac{1}{4}, \frac{1}{2}\right)$ for $s \in B$, and $d(s) = \frac{1}{4}$ for $s \in C$, with $\mu(A), \mu(B), \mu(C) > 0$ and $A \cup B \cup C = \YY$. This is a challenging---yet intriguing---problem; at present, we do not even have a candidate representation for the limiting random variable in this setting.
\end{openquest}

%% file: appendix.tex
\appendix

\section{Some Technical Results}
\label{se:CH3appendix}

In this section, we give some technical results and their proofs. We start with a lemma providing a technical estimate.

\begin{lemma} \label{le:chineq}
The function $c(r,s)$ in \eqref{equality_beta} can be bounded from above as
\begin{equation} \label{eq:inqc(r,s)}
c(r,s)\leq\frac{1}{d(r)}+\frac{1}{1-d(r,s)}
\end{equation}
for $d(r,s)=d(r)+d(s)$ and $d(s) \in (0,\frac{1}{2})$.
\end{lemma}

\begin{proof}
As a function of $d$, the Beta function can be written as
\begin{equation} \label{eq:chtwo}
c(r,s)
=\int_{0}^{\infty}x^{d(r)-1}(x+1)^{d(s)-1}dx
=\int_{0}^{1}x^{d(r)-1}(1-x)^{-d(r,s)}dx.
\end{equation}
Note that $ab \leq a+b$ for $a,b>0$, which implies
\begin{equation} \label{eq:equivalence}
x^{d(r)-1}(1-x)^{-d(r,s)} \leq x^{d(r)-1} + (1-x)^{-d(r,s)}.
\end{equation}
Then, \eqref{eq:equivalence} yields
\begin{equation}
c(r,s) 	
\leq \int_{0}^{1}x^{d(r)-1}dx+\int_{0}^{1}(1-x)^{-d(r,s)}dx
\leq\frac{1}{d(r)}+\frac{1}{1-d(r,s)},
\end{equation}
concluding the proof of the Lemma. \end{proof}

The following lemma ensures that, in the second regime, the conditions of Theorem \ref{th:main_result} imply the desired regularity for the kernels $f$ in \eqref{equality_f^(u,l)_t,(m,s)}.

\begin{lemma} \label{eq:finL2}
Recall the kernel $\mathfrak{f}$ from \eqref{equality_f^(u,l)_t,(m,s)} and let $d(r) \in \left(\frac{1}{4},\frac{1}{2} \right)$. If 
\begin{equation}
    \int_{\mathbb{Y}} \frac{\sigma^2(r)}{(1- 2d(r))(2d(r) - 1/2)} \mu(dr) < \infty,
\end{equation}
then $\mathfrak{f}^{\sigma}(r,s,x,y) = \mathfrak{f}^{(r,s)}(x,y) \sigma(r) \sigma(s)$ takes values in $ L^2(\RR^2:L^2(\YY^2))$.
\end{lemma}

\begin{proof}
First, note that
\begin{align}
&	
\int_{\RR^2} \left| \int_{0}^{1} (v-x_{1})_{+}^{d(r)-1} (v-x_{2})_{+}^{d(s)-1} dv \right|^2 dx_{1}dx_{2} \nonumber	\\
& =
\int_{\RR^2} \left( \int_{0}^{1} (v-x_{1})_{+}^{d(r)-1} (v-x_{2})_{+}^{d(s)-1} dv \right) \left(
\int_{0}^{1} (u-x_{1})_{+}^{d(r)-1} (u-x_{2})_{+}^{d(s)-1} du \right) dx_{1}dx_{2} \nonumber	\\
& =
\int_{0}^{1} \int_{0}^{1} \left(
\int_{\RR} (v-x_{1})_{+}^{d(r)-1} (u-x_{1})_{+}^{d(r)-1} dx_{1} \right) \left(
\int_{\RR} (v-x_{2})_{+}^{d(s)-1} (u-x_{2})_{+}^{d(s)-1} dx_{2} \right) dv du \nonumber	\\
& =
\int_{0}^{1} \int_{0}^{1} \left(
\int_{\RR} (z_{1})_{+}^{d(r)-1} (u-v+z_{1})_{+}^{d(r)-1} dz_{1} \right) \left(
\int_{\RR} (z_{2})_{+}^{d(s)-1} (u-v+z_{2})_{+}^{d(s)-1} dz_{2} \right) dv du,	\label{eq:subz} 
\end{align}
where \eqref{eq:subz} follows by substituting $x_{1}=v-z_{1}$ and $x_{2}=v-z_{2}$. We want to perform the change of variables $z_i = (u - v)  w_i, i = 1,2$. For this, we first see that, for the sets
\begin{equation*}
\begin{split}
    \mathcal{A} &\doteq \{z_1 \ge 0, z_2 \ge 0\} = \{u \ge v, w_1 \ge 0, w_2 \ge 0 \} \cup \{u \le v, w_1 \le 0, w_2 \le 0 \} \\
    \mathcal{B} &\doteq \{u-v + z_1 \ge 0, u -v + z_2 \ge 0\} = \{u \ge v, w_1 \ge -1, w_2 \ge -1 \} \cup \{u \le v, w_1 \le -1, w_2 \le -1 \},
\end{split}
\end{equation*}
we have that
\begin{equation} \label{eq:decomp-a-b}
    \mathcal{A} \cap \mathcal{B} = \{u \ge v, w_1 \ge 0, w_2 \ge 0 \} \cup \{u \le v, w_1 \le -1, w_2 \le -1 \}.
\end{equation}
In view of \eqref{eq:decomp-a-b} and the change of variables $z_i = (u - v)  w_i, i = 1,2$, \eqref{eq:subz} implies that
\begin{align}
&	\int_{\RR^2} \left| \int_{0}^{1} (v-x_{1})_{+}^{d(r)-1} (v-x_{2})_{+}^{d(s)-1} dv \right|^2 dx_{1}dx_{2} \nonumber \\
&=
\int_{0}^{1} \int_{0}^{u} (u-v)^{2d(r,s)-2}  dv du \nonumber\\
&\hspace{2cm}\times \int_0^\infty (w_{1})^{d(r)-1} (1+w_{1})^{d(r)-1} dw_{1} 
\int_0^\infty (w_{2})^{d(s)-1} (1+w_{2})^{d(s)-1} dw_{2}  \label{eq:subuv}	\\
&\quad+\int_{0}^{1} \int_{u}^{1} (v-u)^{2d(r,s)-2}  dv du \nonumber \\
&\hspace{2cm}\times \int_{-\infty}^{-1} (-w_{1})^{d(r)-1} (-1-w_{1})^{d(r)-1} dw_{1} 
\int_{-\infty}^{-1} (-w_{2})^{d(s)-1} (-1-w_{2})^{d(s)-1} dw_{2} \label{eq:subuv2} \\
&=2 \int_0^\infty (w_{1})^{d(r)-1} (1+w_{1})^{d(r)-1} dw_{1} 
\int_0^\infty (w_{2})^{d(s)-1} (1+w_{2})^{d(s)-1} dw_{2} \nonumber\\
&\hspace{2cm} \times \left( \int_0^1 \left( \int_{0}^{u} (u-v)^{2d(r,s)-2} + \int_{u}^{1} (v-u)^{2d(r,s)-2}  \right) dv du    \right) \label{eq:subuv3},
\end{align}
where \eqref{eq:subuv3} follows by the change of variables $x_i = -1 -w_i, i = 1,2$ of the integrals in \eqref{eq:subuv2} (we still denote the new variable by $w_i$). By recalling the definition of $c(r) \doteq c(r,r)$, in \eqref{equality_beta} and elementary calculations, this says that
\begin{equation} \label{eq:subuv4}
\begin{split}
    \int_{\RR^2} \left| \int_{0}^{1} (v-x_{1})_{+}^{d(r)-1} (v-x_{2})_{+}^{d(s)-1} dv \right|^2 dx_{1}dx_{2} &= 4 c(r) c(s) \frac{1}{2d(r,s) - 1} \frac{1}{2d(r,s)} \\
    &\le 4 c(r) c(s) \frac{1}{2d(r,s) - 1},
\end{split}
\end{equation}
since $d(r) > 1/4$.

We see that, by Tonelli's theorem,
\begin{align}
\left\| \mathfrak{f}^\sigma \right\|^2_{L^{2}(\RR^2: {L^2(\YY^2)})}
&=  \int_{\RR^2} \left( \int_{\mathbb{Y}^2} \left| \sigma(r) \sigma(s)  \int_{0}^{1} (v-x_{1})_{+}^{d(r)-1} (v-x_{2})_{+}^{d(s)-1} dv  \right|^2 \mu(dr) \mu(ds)  \right) dx_1 dx_2 \nonumber \\
& \leq 4
\int_{\mathbb{Y}^2} \sigma^2(r) \sigma^2(s) c(r) c(s) \frac{1}{2d(r,s) - 1} \mu(dr) \mu(ds)  \label{al:add1} \\
&\leq 4 \int_{\mathbb{Y}^2} \left(\frac{1}{d(r)} + \frac{1}{1-2d(r)} \right) \left(\frac{1}{d(s)} + \frac{1}{1- 2d(s)} \right) \frac{\sigma^2(r) \sigma^2(s)}{2d(r,s) - 1}	\mu(dr) \mu(ds)  \label{al:add2} \\
& \leq 36  \left( \int_{\mathbb{Y}} \frac{\sigma^2(r)}{(1- 2d(r))(2d(r) - 1/2)} \mu(dr) \right)^2 \label{eq:L3},
\end{align}
where \eqref{al:add1} is due to \eqref{eq:subuv4}, \eqref{al:add2} follows from  \eqref{eq:inqc(r,s)}, and \eqref{eq:L3} follows by recalling that, for $1/4 \le d(r),d(s) \le 1/2$, we have the relations
\[
\frac{1}{2d(r,s) - 1} = \frac{1}{2d(r) - 1/2 + 2d(s) - 1/2} \le \frac{1}{2d(r) - 1/2}\frac{1}{2d(s) - 1/2}, \quad \frac{1}{d(r)} \le  \frac{2}{1-2d(r)}.
\]
The proof is concluded upon noticing that the quantity in \eqref{eq:L3} is finite.
\end{proof}

The following key lemma provides the necessary ingredients to obtain convergence to double Wiener-It\^o integrals with sample paths in $L^2(\YY^2)$. It generalizes Proposition 14.3.2 in \citet{giraitis}.
\begin{lemma} \label{Lemma_off_diagonals}
Let $C_N$ be as in \eqref{eq:CN-vector}. Consider a linear combination of an off-diagonal tuple
\begin{equation*}
\quad Q_{2}(C_{N}) \in L^2(\mathbb{Y}^2:\RR^{(H+1)}), \quad Q_{2}(C_{N})(r,s) \doteq \sum_{j_{1} \neq j_{2}} C_{N}(j_{1},j_{2},r,s)\varepsilon_{j_{1}}(r) \varepsilon_{j_{2}}(s) .
\end{equation*}
Assume that there exists a kernel $f \in \mathcal{H}^2$, such that, denoting $f^\sigma(r,s,x,y) \doteq f(r,s,x,y) \sigma(r) \sigma(s)$ and $\bar f \doteq (f,\dots, f) \in L^{2}(\RR^{2}: L^2(\mathbb{Y}^2:\RR^{(H+1)}))$, the functions defined in \eqref{eq:def-cntilde} satisfy, as $N \to \infty$,
\begin{equation} \label{eq:lemma-gira-ass}
\| \widetilde{C}_{N}^\sigma-\bar f^\sigma \|_{L^{2}(\RR^{2}: L^2(\mathbb{Y}^2:\RR^{(H+1)}))} \to 0.
\end{equation}
Then,
\begin{equation*}
    Q_{2}(C_{N}) \overset{\distr}{\to} \mathcal{I}_{2}(\bar f) \doteq \begin{pmatrix}
        \mathcal{I}_{2}^\sigma(f) \\
        \vdots \\
        \mathcal{I}_{2}^\sigma(f)
    \end{pmatrix} , \quad \text{weakly in } L^2(\mathbb{Y}^2:\RR^{(H+1)}),
\end{equation*}
where each coordinate of the limit is understood in the sense of Definition \ref{def:double-in-banach}.
\end{lemma}

\begin{proof}
For notational simplicity, we write from here on $\mathcal{I}_{2}(f)$ instead of $\mathcal{I}_{2}^\sigma(f)$.
By a truncation argument, it is enough to prove that for all $\varepsilon > 0$, there exists a special kernel $f_{\varepsilon}$ (cf. Section \ref{se:biW}), its vectorization $\bar f_{\varepsilon}$, and a corresponding $C_{N,\varepsilon}$, such that
\begin{align}
\Var \|Q_{2}(C_{N})-Q_{2}(C_{N,\varepsilon}) \|_{L^2(\mathbb{Y}^2:\RR^{(H+1)})} \leq \varepsilon,	\label{Con1}		\\
\Var \| \mathcal{I}_{2}(\bar f_{\varepsilon})- \mathcal{I}_{2}(\bar f) \|_{L^2(\mathbb{Y}^2:\RR^{(H+1)})} \leq \varepsilon, 	\label{Con2}	\\	
Q_{2}(C_{N,\varepsilon}) \overset{d}{\longrightarrow} \mathcal{I}_2(\bar f_{\varepsilon}) ,	\label{Con3}
\end{align}
as $N \to \infty$, where the weak convergence in \eqref{Con3} is in the topology of $L^2(\mathbb{Y}^2:\RR^{(H+1)})$ and $C_{N,\varepsilon}$ is defined by
\begin{equation} \label{eq:CNE-def}
C_{N,\varepsilon}(j_{1},j_{2},r,s) 
\doteq 
N^{-1} \bar f_{\varepsilon}^{(r,s)} \left( \frac{j_{1}}{N},\frac{j_{2}}{N}\right), \hspace*{0.2cm} r,s \in \mathbb{Y}.
\end{equation}
Indeed, if \eqref{Con1}--\eqref{Con3} are true, then we can write, for each $\varepsilon > 0$
\begin{equation}
\begin{aligned}
    Q_2(C_N) - \mathcal{I}_2(\bar f) &= Q_2(C_N) - Q_2(C_{N,\varepsilon}) + Q_2(C_{N,\varepsilon}) -  \mathcal{I}_2(\bar f_\varepsilon) +  \mathcal{I}_2(\bar f_\varepsilon) - \mathcal{I}_2(\bar f)  \\
    &= \mathcal{A}_{N,\varepsilon} + \mathcal{B}_{N,\varepsilon} + \mathcal{C}_{N,\varepsilon},
\end{aligned}
\end{equation}
where
\begin{equation}
    \mathcal{A}_{N,\varepsilon} \doteq Q_2(C_N) - Q_2(C_{N,\varepsilon}), \quad \mathcal{B}_{N,\varepsilon} \doteq Q_2(C_{N,\varepsilon}) - \mathcal{I}_2(\bar f_\varepsilon), \quad \mathcal{C}_{N,\varepsilon} \doteq \mathcal{I}_2(\bar f_\varepsilon) - \mathcal{I}_2(\bar f) .
\end{equation}
Then, by Condition \eqref{Con1} it follows that $\mathcal{A}_{N,\varepsilon}  \xrightarrow{P} 0$ in $L^2(\mathbb{Y}^2:\RR^{(H+1)})$. Moreover, by Condition \eqref{Con3}, it follows that $\mathcal{B_{N,\varepsilon}} \to 0$ in distribution, and so in probability. Finally, \eqref{Con2} implies that $\mathcal{C_{N,\varepsilon}} \xrightarrow{P} 0$. We can then conclude that
\begin{equation}
    Q_2(C_N) \to \mathcal{I}_2(\bar f), \quad \text{weakly in } L^2(\mathbb{Y}^2:\RR^{(H+1)}),
\end{equation}
from Slutsky's lemma, by sending first $N \to \infty$, and then $\varepsilon \to 0$.

We turn to verifying conditions \eqref{Con1}--\eqref{Con3}. First note that, whenever $i_1 \neq i_2, j_1 \neq j_2$,
\begin{equation*}
\EE (\varepsilon_{i_{1}}(r_{1}) \varepsilon_{j_{1}}(s_{1}) \varepsilon_{i_{2}}(r_{2}) \varepsilon_{j_{2}}(s_{2}))=
\begin{cases}
\sigma(r_{1},s_{1}) \sigma(r_{2},s_{2}),			&\text{ if } i_{1}=j_{1}, i_{2}=j_{2}, \\
\sigma(r_{1},s_{2}) \sigma(s_{1},r_{2}),			&\text{ if } i_{1}=j_{2}, j_{1}=i_{2}, \\
0, 										&\text{ otherwise.}
\end{cases}
\end{equation*}
Then, for \eqref{Con1},
\begin{align}
&\EE \|Q_{2}(C_{N})\|^2_{L^2(\mathbb{Y}^2:\RR^{(H+1)})} \nonumber
\\ & =	
\EE \Big(\int_{\mathbb{Y}} \int_{\mathbb{Y}} 
\bigg\|\sum_{j_{1} \neq j_{2}} C_{N}(j_{1},j_{2},r,s)\varepsilon_{j_{1}}(r)\varepsilon_{j_{2}}(s)\bigg\|_{\RR^{H+1}}^2 \mu(dr) \mu(ds)	\Big)	\nonumber
\\ & =	
\int_{\mathbb{Y}} \int_{\mathbb{Y}} \sum_{j_{1} \neq j_{2}} \sum_{i_{1} \neq i_{2}} (C_{N}(j_{1},j_{2},r,s))' C_{N}(i_{1},i_{2},r,s) \EE(\varepsilon_{j_{1}}(r) \varepsilon_{j_{2}}(s) \varepsilon_{i_{1}}(r) \varepsilon_{i_{2}}(s)) \mu(dr) \mu(ds)	\nonumber
\\ & =	
\int_{\mathbb{Y}} \int_{\mathbb{Y}} \sum_{j_{1} \neq j_{2}} \| C_{N}(j_{1},j_{2},r,s) \|^2 \nonumber
 \sigma^2(r)\sigma^2(s) \mu(dr) \mu(ds) \\
&\quad \quad + \int_{\mathbb{Y}} \int_{\mathbb{Y}} \sum_{j_{1} \neq j_{2}} | C_{N}(j_{1},j_{2},r,s)|' | C_{N}(j_{2},j_{1},r,s)|  \sigma^2(r,s)
 \mu(dr) \mu(ds) \nonumber
\\ & \leq	
2 \int_{\mathbb{Y}} \int_{\mathbb{Y}} \int_{\RR^2}  N^2  \| C_{N}(\lceil x_{1}N \rceil ,\lceil x_{2}N \rceil ,r,s) \|^2 dx_{1} dx_{2} \sigma^2(r)\sigma^2(s) \mu(dr) \mu(ds)	\label{eq:890}
\\ & \le		
2 \| \widetilde{C}^{\sigma}_{N} \|^2_{L^{2}(\RR^2: L^2(\mathbb{Y}^2:\RR^{(H+1)}))} \label{eq:890+1}
\end{align}
with $\widetilde{C}^{\sigma}_{N}(x_1,x_2,r,s)$ as in \eqref{eq:def-cntilde}, and 
where \eqref{eq:890} holds by noticing that $\sigma^2(r,s) \le \sigma^2(r) \sigma^2(s)$.
This implies
\begin{align}
\EE\|Q_{2}(C_{N})-Q_{2}(C_{N,\varepsilon})\|^2_{L^{2}(\RR^2: L^2(\mathbb{Y}^2:\RR^{(H+1)}))} 
&= 
\EE\|Q_{2}(C_{N} - C_{N,\varepsilon})\|^2_{L^{2}(\RR^2: L^2(\mathbb{Y}^2:\RR^{(H+1)}))} \nonumber \\
& \leq			
2 \| \widetilde{C}^{\sigma}_{N} - \widetilde{C}^{\sigma}_{N,\varepsilon} \|^2_{L^{2}(\RR^2: L^2(\mathbb{Y}^2:\RR^{(H+1)}))},
\label{eq:890+2}
\end{align}
where 
\begin{equation}
    \widetilde{C}^{\sigma}_{N,\varepsilon}(x_{1},x_{2},r,s) \doteq N C_{N,\varepsilon}(x_1,x_2,r,s) \sigma(r) \sigma(s).
\end{equation}
It remains to bound the right hand side of \eqref{eq:890+2}.
First note that if $f \in \mathcal{H}^2$ is special, then so is $f^\sigma$. Then, for any $N_0 \ge 1$ and any special $\bar f^\sigma_\varepsilon, \widetilde C^\sigma_{N,\varepsilon}$, by the triangle and Cauchy-Schwarz inequalities,
\begin{align}
&
\| \widetilde{C}^{\sigma}_{N} -\widetilde{C}^{\sigma}_{N,\varepsilon} \|^2_{L^{2}(\RR^2: L^2(\mathbb{Y}^2:\RR^{(H+1)}) )}
\leq
3\| \widetilde{C}^{\sigma}_{N} -\widetilde{C}^{\sigma}_{N_{0}}\|^2_{L^{2}(\RR^2: L^2(\mathbb{Y}^2:\RR^{(H+1)}))} \nonumber
\\& \hspace{2cm}+
3\| \widetilde{C}_{N_{0}}^{\sigma} -\bar f^{\sigma}_{\varepsilon}\|^2_{L^{2}(\RR^2: L^2(\mathbb{Y}^2:\RR^{(H+1)}) )}+
3\| \bar f^{\sigma}_{\varepsilon} -\widetilde{C}^{\sigma}_{N,\varepsilon}\|^2_{L^{2}(\RR^2: L^2(\mathbb{Y}^2:\RR^{(H+1)}))}. \label{eq:301}
\end{align}
By assumption, there is a $N_{0} \geq 1$ such that, for the choice of kernels $f^\sigma$
\begin{align} 
&
\| \widetilde{C}^{\sigma}_{N} -\widetilde{C}^{\sigma}_{N_{0}}\|^2_{L^{2}(\RR^2: L^2(\mathbb{Y}^2:\RR^{(H+1)}))} \nonumber
\\&\leq			
2\| \widetilde{C}^{\sigma}_{N} -\bar f^{\sigma} \|^2_{L^{2}(\RR^2: (L^2(\mathbb{Y}^2)^{\times (H+1)}))}+
2\| \bar f^{\sigma}-\widetilde{C}^{\sigma}_{N_{0}}\|^2_{L^{2}(\RR^2: L^2(\mathbb{Y}^2:\RR^{(H+1)}))}
\leq
\frac{\varepsilon}{18} \label{eq:302}
\end{align}
for all $N \geq N_{0}$.
Given $N_{0}\geq 1$ and $\varepsilon>0$, since $f \in \mathcal{H}^2$, there exist special kernels $f_{\varepsilon}\in \mathcal{H}^2$ as in \eqref{eq:simple-kernel} (with $c$ replaced by $\bar c$) and their induced vectorization $\bar f_{\varepsilon}^\sigma$ from Theorem \ref{lemma:simple-dominated-above} such that
\begin{equation} \label{eq:303}
\| \widetilde{C}^{\sigma}_{N_{0}} -\bar f^{\sigma}_{\varepsilon} \|^2_{L^{2}(\RR^2: L^2(\mathbb{Y}^2:\RR^{(H+1)}) )} 
\leq 
\| \widetilde{C}^\sigma_{N_{0}} -\bar f^{\sigma} \|^2_{L^{2}(\RR^2: L^2(\mathbb{Y}^2:\RR^{(H+1)}))} + \| \bar f^{\sigma} - \bar f^{\sigma}_{\varepsilon} \|^2_{L^{2}(\RR^2: L^2(\mathbb{Y}^2:\RR^{(H+1)}))} \le \varepsilon/9.
\end{equation}
Moreover, the function $\widetilde{C}_{N,\varepsilon}$ derived from $C_{N,\varepsilon}$ satisfies, for all $\varepsilon > 0$,
\begin{align*} 
&
\|\bar f^{\sigma}_{\varepsilon}-\widetilde{C}^{\sigma}_{N,\varepsilon}\|^2_{L^{2}(\RR^2: L^2(\mathbb{Y}^2:\RR^{(H+1)}))}
\\&=		
\int_{\RR^2} \left\| \bar f^{\sigma}_{\varepsilon}(x_{1},x_{2})-\bar f^{\sigma}_{\varepsilon} \Big(\frac{\lceil  x_{1}N \rceil }{N} , \frac{\lceil x_{2}N \rceil }{N} \Big) \right\|^2_{L^2(\mathbb{Y}^2:\RR^{(H+1)})} dx_{1}dx_{2}
\\&\le (H+1)		
\int_{\YY^2} \sigma^2(r) \sigma^2(s) \| \bar c_{i_1,i_2}(r,s)\|^2 \mu(dr) \mu(ds)  \\
&\quad\quad \times
\int_{\RR^2} \left( \mathds{1}_{\Delta_{i_1}}(x_{1}) \mathds{1}_{\Delta_{i_2}}(x_{2}) - \mathds{1}_{\Delta_{i_1}}\Big(\frac{\lceil x_{1}N \rceil }{N} \Big)
\mathds{1}_{\Delta_{i_2}}\Big(\frac{\lceil x_{2}N \rceil }{N} \Big) \right) dx_{1}dx_{2} \to  0,
\end{align*}
as $N \to \infty$, where we used that $f_{\varepsilon}^{\sigma}$ is a special kernel and DCT. Therefore, we can choose an $N_1 \ge N_0$ such that, for all $N \ge N_1$,
\begin{equation} \label{eq:304}
\|\bar f^{\sigma}_\varepsilon -\widetilde{C}^{\sigma}_{N,\varepsilon} \|^2_{L^{2}(\RR^2: L^2(\mathbb{Y}^2:\RR^{(H+1)}) )}\\
\leq
\varepsilon /9.
\end{equation}
Combining the relations in \eqref{eq:301}, \eqref{eq:302}, \eqref{eq:303}, and \eqref{eq:304} finishes the proof of
\eqref{Con1}. 

We move to Condition \eqref{Con2}. Note that
\begin{align}
&
\Var\| \mathcal{I}_2(\bar f_{\varepsilon})- \mathcal{I}_2(\bar f) \|_{L^2(\mathbb{Y}^2:\RR^{(H+1)})}	\nonumber
\\&=
\EE \int_{\YY^2} \int_{\RR^4} (\bar f^{(r,s)}(x_1,x_2)- \bar f_{\varepsilon}^{(r,s)}(x_1,x_2))' ( \bar f^{(r,s)}(y_1,y_2)- \bar f_{\varepsilon}^{(r,s)}(y_1,y_2))  \nonumber \\
&\quad \quad \quad  \quad \quad \quad \times W^{(r)}(dx_1) W^{(s)}(dx_2) W^{(r)}(dy_1) W^{(s)}(dy_2) \mu(dr) \mu(ds) \nonumber
\\&= \int_{\YY^2} \int_{\RR^2} \| \bar f^{(r,s)}(x_1,x_2)- \bar f_{\varepsilon}^{(r,s)}(x_1,x_2) \|^2 dx_1 dx_2 \sigma^2(r) \sigma^2(s) \mu(dr) \mu(ds) \nonumber \\
& \quad \quad + \int_{\YY^2} \int_{\RR^2} (\bar f^{(r,s)}(x_1,x_2)- \bar f_{\varepsilon}^{(r,s)}(x_1,x_2))' (\bar f^{(r,s)}(x_2,x_1)- \bar f_{\varepsilon}^{(r,s)}(x_2,x_1)) dx_1 dx_2 \nonumber
\\& \quad \quad \quad \quad \times \sigma^2(r,s)  \mu(dr) \mu(ds) \nonumber \\
&\le 2\int_{\YY^2} \int_{\RR^2} \| \bar f^{(r,s)}(x_1,x_2)- \bar f_{\varepsilon}^{(r,s)}(x_1,x_2) \|^2 dx_1 dx_2 \sigma^2(r) \sigma^2(s) \mu(dr) \mu(ds),
\label{al:two_sep}
\end{align}
where the last inequality follows by Young's inequality. Moreover, by Theorem \ref{lemma:simple-dominated-above} we can choose $f_\varepsilon^{(r,s)}(x,y)$ to be dominated from $f^{(r,s)}(x,y)$ pointwise. This says that, from Lemma \ref{eq:finL2} and Cauchy-Schwarz,
\begin{equation*}
    \Var\| \mathcal{I}_2(\bar f_{\varepsilon})- \mathcal{I}_2(\bar f) \|_{L^2(\mathbb{Y}^2:\RR^{(H+1)})} \le 
    3 \| \bar f^{\sigma}- \bar f_{\varepsilon}^{\sigma} \|^2_{L^{2}(\RR^2: L^2(\mathbb{Y}^2:\RR^{(H+1)}) )}.
\end{equation*}
Then, by DCT, we can select $f_{\varepsilon}$ such that, in addition to \eqref{Con1},
\begin{equation}
    \Var\| \mathcal{I}_2(\bar f_{\varepsilon})-\mathcal{I}_2(\bar f) \|_{L^2(\mathbb{Y}^2:\RR^{(H+1)})}	 \le \varepsilon.
\end{equation}
For \eqref{Con3}, note that
\begin{align*}
Q_{2}(C_{N,\varepsilon})
& =
\sum_{j_{1} \neq j_{2}} C_{N, \varepsilon}(j_{1},j_{2},\cdot,\cdot)\varepsilon_{j_{1}}(\cdot) \varepsilon_{j_{2}}(\cdot)  =
\sum_{j_{1} \neq j_{2}} N^{-1} \bar f_{\varepsilon}^{(\cdot,\cdot)} \Big(\frac{j_{1}}{N},\frac{j_{2}}{N}\Big) \varepsilon_{j_{1}}(\cdot) \varepsilon_{j_{2}}(\cdot).
\end{align*}
In order to show that $Q_{2}(C_{N,\varepsilon})$ converges weakly to $\mathcal{I}_2(\bar f_\varepsilon)$ in $L^2(\mathbb{Y}^2:\RR^{(H+1)})$, we employ Theorem 2 in \cite{CreKad86}. We rephrase the conditions for weak convergence here. A sequence $\xi_n$  converges weakly to $\xi \in L^2(\mathbb{Y}^2:\RR^{(H+1)})$ if
\begin{enumerate}[label=\textit{(\roman*)}]
    \item the finite-dimensional distributions converge, i.e.,
    \begin{equation}
        (\xi_n(r_1,s_1), \dots, \xi_n(r_q, s_q))' \overset{d}{\to} 
        (\xi(r_1,s_1), \dots, \xi(r_q, s_q))',
    \end{equation}
    in $\RR^{q\times (H+1)}$ for all $r_1,\dots,r_q,s_1,\dots,s_q \in \YY $ and $q \in \NN$.
    \label{item:kadelka1}
    \item For all $r,s \in \YY$ we have that $\EE \| \xi_n(r,s) \|^2_{\RR^{(H+1)}} \to \EE \| \xi(r,s) \|^2_{\RR^{(H+1)}}$.
    \label{item:kadelka2}
    \item there exists a $\mu \otimes \mu$-integrable function $g: \YY^2 \to \RR$ such that $ \EE \|\xi_{n}(r,s) \|_{\RR^{H+1}}^2 \leq g(r,s)$ for each $ (r,s) \in \YY^2$.
    \label{item:kadelka3}
\end{enumerate}

\textit{Proof of condition \ref{item:kadelka1}:}
Since $\bar f_{\varepsilon}$ are special, they take non-zero values in a finite number of intervals $\{\Delta_i\}_{i=1,\dots,J}$, and so can be rewritten as
\begin{equation}
    \bar f_{\varepsilon}^{(r,s)}(x,y) = \sum_{\substack{i_1,i_2 = 1 \\ i_1 \neq i_2}}^J \bar c_{i_1,i_2}(r,s) \mathds{1}_{\Delta_{i_1}}(x) \mathds{1}_{\Delta_{i_2}}(y),
\end{equation}
with suitable functions $c_{i_1,i_2} : \YY^2 \to \RR$, $\bar c = (c,\dots,c) \in L^2 (\YY^2, \RR^{H+1})$. Then, for $r,s \in \YY$, from \eqref{eq:CNE-def},
\begin{equation} \label{eq:q2cn-cont}
\begin{split}
Q_{2}(C_{N,\varepsilon})(r,s) &= 
\sum_{\substack{i_1,i_2 = 1 \\ i_1 \neq i_2}}^J \bar c_{i_1, i_2}(r,s) N^{-1} 
\sum_{j_{1} \neq j_{2}}\varepsilon_{j_{1}}(r) \varepsilon_{j_{2}}(s)
\mathds{1} _{ \{ \frac{j_{1}}{N} \in \Delta_{1}, \frac{j_{2}}{N} \in \Delta_{2} \}} \\
& =
\sum_{\substack{i_1,i_2 = 1 \\ i_1 \neq i_2}}^J \bar c_{ i_1, i_2}(r,s)
W^{(r)}_{N} (\Delta_{i_1}) W^{(s)}_{N} (\Delta_{i_2}),
\end{split}
\end{equation}
where
\begin{equation*}
W^{(r)}_{N}(\Delta_{i}) \doteq N^{-\frac{1}{2}} \sum_{j : \frac{j}{N}\in \Delta_{i}} \varepsilon_{j}(r) 
= 
N^{-\frac{1}{2}} \sum_{j=1}^N \varepsilon_j(r) \mathds{1}_{\{\frac{j}{N} \in \Delta_i\}}.
\end{equation*}
Similarly, 
\begin{equation}
    \begin{pmatrix}
        Q_{2}(C_{N,\varepsilon})(r_1,s_1) \\
        \vdots \\
        Q_{2}(C_{N,\varepsilon})(r_q,s_q)
    \end{pmatrix} = \begin{pmatrix}
        \sum_{i_1 \neq i_2}^J \bar c_{ i_1, i_2}(r_1,s_1)
W^{(r_1)}_{N} (\Delta_{i_1}) W^{(s_1)}_{N} (\Delta_{i_2}) \\
\vdots \\
\sum_{i_1 \neq i_2}^J \bar c_{ i_1, i_2}(r_q,s_q)
W^{(r_q)}_{N} (\Delta_{i_1}) W^{(s_q)}_{N} (\Delta_{i_2})
    \end{pmatrix}.
\end{equation}
We study the joint weak convergence of $W_{N}^{(x)}(\Delta_{p})$, for $p=1,\dots, J$ and $x = r_1,\dots,r_q,s_1,\dots,s_q$. Denote
\begin{equation} \label{eq:WN-def}
    W_N \doteq \begin{pmatrix}
        W_N(\Delta_1) \\
        \vdots \\
        W_N(\Delta_J)
    \end{pmatrix}, 
    \quad W_N(\Delta_p) \doteq 
    \begin{pmatrix}
        W^{(r_1)}_{N}(\Delta_p) \\ W^{(s_1)}_{N}(\Delta_p) \\
        \vdots   \\
        W^{(r_q)}_{N}(\Delta_p) \\ W^{(s_q)}_{N}(\Delta_p)
    \end{pmatrix}, 
    \quad p = 1,\dots J.
\end{equation}
Note that from a standard multivariate CLT, we have, for $p=1,\dots, J$,
\begin{equation*}
    W_N(\Delta_p) \xrightarrow{d} W(\Delta_p) \in \RR^{2q}, \quad \Sigma_p  \doteq \begin{pmatrix}
        \sigma^2(r_1) & \sigma(r_1,s_1) & \sigma(r_1,r_2) & \hdots & \sigma(r_1,s_q) \\
        \vdots & \vdots & \vdots & \ddots & \vdots \\
        \sigma(s_q,r_1) & \sigma(s_q,s_1) & \sigma(s_q,r_2) & \hdots &\sigma^2(s_q)
    \end{pmatrix},
\end{equation*}
where $W(\Delta_p)$ is a Gaussian random vector with covariance matrix $\Sigma_p$. Moreover, $W(\Delta_p)$ and $W(\Delta_s)$ are independent whenever $p \neq s$ since $\Delta_p$ and $\Delta_s$ are disjoint. We therefore obtain the joint convergence, as $N \to \infty$,
\begin{equation}
    W_N \xrightarrow{d} W
\end{equation}
where $W$ is a $\RR^{J \times (2q)}$-valued Gaussian random variable that has a covariance matrix given by the block diagonal form $\Sigma_{m,n} = \Sigma_m \delta_{m,n}, m , n= 1,\dots J$.

Then, we can recast
\begin{equation}
    \begin{pmatrix}
        Q_{2}(C_{N,\varepsilon})(r_1,s_1) \\
        \vdots \\
        Q_{2}(C_{N,\varepsilon})(r_q,s_q)
    \end{pmatrix} = D (W_N), 
\end{equation}
where $W_N$ is as in \eqref{eq:WN-def} and $D$ is a suitable function involving $c_{i_1,i_2}, i_1, i_2 = 1, \dots J$ (see, \cite{giraitis} pp. 535-536). We then obtain by the continuous mapping theorem that, as $N \to \infty$,
\begin{equation*}
    \begin{pmatrix}
        Q_{2}(C_{N,\varepsilon})(r_1,s_1) \\
        \vdots \\
        Q_{2}(C_{N,\varepsilon})(r_q,s_q)
    \end{pmatrix} \xrightarrow{d} \begin{pmatrix}
        \sum_{i_1 \neq i_2} \bar c_{i_1,i_2}(r_1,s_1)
W^{(r_1)} (\Delta_{i_1}) W^{(s_1)} (\Delta_{i_2}) \\
\vdots \\
\sum_{i_1 \neq i_2} \bar c_{i_1,i_2}(r_q,s_q)
W^{(r_q)} (\Delta_{i_1}) W^{(s_q)} (\Delta_{i_2})
    \end{pmatrix} = \begin{pmatrix}
    \mathcal{I}_2(\bar f_\varepsilon)(r_1,s_1) \\
    \vdots \\
    \mathcal{I}_2(\bar f_\varepsilon)(r_q,s_q)
\end{pmatrix} \in \RR^{q \times (H+1)}.
\end{equation*}

\textit{Proof of condition \ref{item:kadelka2}:}
For fixed $r,s \in \YY$, we have
\begin{align*}
    &\EE \|\mathcal{I}_2(\bar f_{\varepsilon})(r,s) \|_{\RR^{(H+1)}}^2 
    \\&= 
    \EE \left( \int_{\RR^4} ( \bar f^{(r,s)}_{\varepsilon}(x_1,x_2) )' \bar f^{(r,s)}_{\varepsilon}(y_1,y_2) W^{(r)}(dx_1) W^{(s)}(dx_2) W^{(r)}(dy_1) W^{(s)}(dy_2) \right) \\
    &= \int_{\RR^4} (\bar f^{(r,s)}_{\varepsilon}(x_1,x_2))' \bar f^{(r,s)}_{\varepsilon}(y_1,y_2) \EE \left( W^{(r)}(dx_1) W^{(s)}(dx_2) W^{(r)}(dy_1) W^{(s)}(dy_2) \right) \\
    &= \int_{\RR^2} (\bar f^{(r,s)}_{\varepsilon}(x_1,x_2))' \bar f^{(r,s)}_{\varepsilon}(x_1,x_2)  \EE \left( W^{(r)}(dx_1)  W^{(r)}(dx_1) \right) \EE \left( W^{(s)}(dx_2) W^{(s)}(dx_2) \right) \\
    &\quad \quad + \int_{\RR^2} (\bar f^{(r,s)}_{\varepsilon}(x_1,x_2) )'\bar f^{(r,s)}_{\varepsilon}(x_2,x_1) \EE \left( W^{(r)}(dx_1)  W^{(s)}(dx_1) \right) \EE \left( W^{(r)}(dx_2) W^{(s)}(dx_2) \right) \\
    &= \int_{\RR^2} (\bar f^{(r,s)}_{\varepsilon}(x_1,x_2))' \bar f^{(r,s)}_{\varepsilon}(x_1,x_2) dx_1 dx_2 \sigma^2(r) \sigma^2(s) \\
    &\quad \quad + \int_{\RR^2} (\bar f^{(r,s)}_{\varepsilon}(x_1,x_2) )' \bar f_{\varepsilon}^{(r,s)}(x_2,x_1) dx_1 dx_2 (\sigma(r,s))^2,
\end{align*}
since $f^{(r,s)}_{\varepsilon}(x,x) = 0$ for all $x$ and from \eqref{eq:covWW}. Then, calculations similar to \eqref{eq:890+1} show that
\begin{multline}
    \EE \| Q_{2}(C_{N,\varepsilon})(r,s) \|_{\RR^{H+1}}^2
    = 
    \int_{\RR^2} \| \widetilde{C}_{N,\varepsilon}(x_1,x_2,r,s) \|^2 \sigma^2(r) \sigma^2(s) dx_1 dx_2 
    \\+
    \int_{\RR^2} (\widetilde{C}_{N,\varepsilon}(x_1,x_2,r,s) )' \widetilde{C}_{N,\varepsilon}(x_2,x_1,r,s) \sigma^2(r,s) dx_1 dx_2
    \to \EE \|\mathcal{I}_2(
    \bar f_{\varepsilon})(r,s) \|^2_{\RR^{H+1}} ,
\end{multline}
where the convergence holds from calculations similar to the ones leading to \eqref{eq:304} and DCT.

\textit{Proof of condition \ref{item:kadelka3}:} Identical calculations as in \eqref{eq:890} show that
\begin{equation} \label{eq:909}
    \EE \| Q_2(C_{N,\varepsilon})(r,s) \|_{\RR^{H+1}}^2
    \le 
    2 \| \widetilde C^{\sigma}_{N,\varepsilon} (r,s) \|^2_{ L^2(\RR^2:\RR^{(H+1)})}.
\end{equation}
Moreover, we have that, from \eqref{eq:303}-\eqref{eq:304}
\begin{equation} \label{eq:910}
\begin{split}
    &
    \| \widetilde C^\sigma_{N,\varepsilon} \|^2_{L^2(\RR^2: L^2(\mathbb{Y}^2:\RR^{(H+1)}))} 
    \\&\le 
    8 \left( \| \widetilde C^\sigma_{N,\varepsilon} - \bar f^\sigma_{\varepsilon} \|^2_{L^2(\RR^2: L^2(\mathbb{Y}^2:\RR^{(H+1)}))} + 
    \| \bar f^\sigma_{\varepsilon} - \bar f^\sigma \|^2_{L^2(\RR^2: L^2(\mathbb{Y}^2:\RR^{(H+1)}))}   + \| \bar f^\sigma \|^2_{L^2(\RR^2: L^2(\mathbb{Y}^2:\RR^{(H+1)}))}\right) \\
    &< \varepsilon + \| \bar f^\sigma \|^2_{L^2(\RR^2: L^2(\mathbb{Y}^2:\RR^{(H+1)}))},
\end{split}
\end{equation}
and so $C^{\sigma}_{N,\varepsilon}(r,s)$ is $\mu \otimes \mu$- integrable.
The proof of \ref{item:kadelka3} follows from \eqref{eq:909}, \eqref{eq:910}, and Lemma \ref{eq:finL2}. It also completes the proof of the Lemma.
\end{proof}